\documentclass{article}

\usepackage{amsfonts,amssymb}
\usepackage{amsmath,bbm}
\usepackage{theorem}          
\usepackage{latexsym}
\usepackage{float}                   % steuert Gleitobjekte wie Figuren,
\usepackage{times}
\usepackage{cite}

\newtheorem{tht}{Theorem}[section]

\newtheorem{thl}[tht]{Lemma}
\newtheorem{thp}[tht]{Proposition}
\newtheorem{thc}[tht]{Corollary}

\newcommand{\ang}{\raisebox{0.2ex}{\scriptsize$\triangleright$}}

\newcommand{\mn}{\medskip\noindent}

\newcommand{\sn}{\smallskip\noindent}

\newcommand{\D}{{\mathcal{D}}} 
 
\newcommand{\cR}{{\mathcal{R}}} 
 
\newcommand{\Hh}{{\mathcal{H}}} 
\newcommand{\cO}{{\mathcal{O}}}

\newcommand{\K}{{\mathcal{K}}}
\newcommand{\U}{{\mathcal{U}}}   \newcommand{\cU}{{\mathcal{U}}}

\newcommand{\X}{{\mathcal{X}}}

\newcommand{\F}{{\mathcal{F}}}

\newcommand{\cS}{\mathcal{S}}

\newcommand{\fA}{\mathfrak{A}}

\newcommand{\gM}{{\mathfrak{M}}}
\newcommand{\g}{{\mathfrak{g}}}

\newcommand{\Z}{\mathbb{Z}}
\newcommand{\N}{\mathbb{N}}
\newcommand{\R}{\mathbb{R}}
\newcommand{\C}{\mathbb{C}}
\newcommand{\FF}{\mathbb{F}}
\newcommand{\FFD}{\mathbb{F}(D)}
\newcommand{\FFDA}{\mathbb{F}(D_\fA)}
\newcommand{\BB}{\mathbb{B}}
\newcommand{\BBA}{\mathbb{B}_1(\fA)}
\newcommand{\FS}{\F(\sigma(\bar y))}
\newcommand{\Sy}{\cS(\sigma(\bar y))}
\newcommand{\dy}{\delta_{k}(\bar y)}
\newcommand{\py}[1][n]{\psi_{{#1}}(\bar y)}

\newcommand{\im}{\mathrm{i}}
\newcommand{\Lin}{{\mathrm{Lin}}}

\newcommand{\ad}{{\mathrm{ad}}}      %adjungierte Wirkung

\newcommand{\Ker}{\mathrm{ker}\,}

\newcommand{\vare}{\varepsilon}

\newcommand{\ov}{\overline}
\newcommand{\rf}[1][]{\textup{\eqref{#1}}}

\newcommand{\sut}{{\cU}_q(\mathrm{su}_{1,1})}       %su(1,1)

\newcommand{\sltn}{{\cU}_q(\mathrm{sl}_{n+1})}
\newcommand{\sutn}{{\cU}_q(\mathrm{su}_{n,1})} 
\newcommand{\sutm}{{\cU}_q(\mathrm{su}_{n})} 
\newcommand{\tr}[1][]{\mathrm{Tr}_{{#1}}\,}

\newcommand{\ip}[2]{\langle {#1},{#2}\rangle}
\newcommand{\qd}{\cO_q(\mathrm{U})}                    %quantum disc
\newcommand{\qm}{\cO_q(\mathrm{Mat}_{n,1})}

\newcommand{\ld}{\mathcal{L}^+(D)}            %L+(D)
\newcommand{\ldop}[1][D_\fA]{\mathcal{L}^+({#1})} %\ldop ODER: \ldop[xxx]
\newcommand{\ldd}{\mathcal{L}(D,D^+)}
\newcommand{\ldda}{\mathcal{L}(D_\fA,D_\fA^+)}
\newcommand{\lddp}{\mathcal{L}(D,D^\prime)}
\newcommand{\dd}{\mathrm{d}}

\newcommand{\Dq}{\mathrm{D}_q}
\newcommand{\Dqq}{\mathrm{D}_{q^2}}
\newcommand{\ra}{r(\bar{T})}
\newcommand{\Aj}{A_{j\pm1}}
\newcommand{\rj}{\rho_{j\pm1}}

\title{Invariant integration theory on non-compact quantum spaces: 
       Quantum $(n,1)$-matrix ball }

\author{Klaus-Detlef K\"ursten  and Elmar Wagner\\
\small{ Fakult\"at f\"ur Mathematik und Informatik}\\ 
\small{ Universit\"at Leipzig, Augustusplatz 10, 04109 Leipzig, Germany}\\ 
\small{E-mail: kuersten@mathematik.uni-leipzig.de / 
wagner@mathematik.uni-leipzig.de} }
\begin{document}
\date{}
\maketitle
\renewcommand{\theenumi}{\roman{enumi}}
\begin{abstract}
An operator theoretic approach to 
invariant integration theory on non-compact quantum spaces 
is introduced on the example of the quantum $(n,1)$-matrix ball $\qm$. 
In order to prove the existence of an 
invariant integral, operator algebras are associated to $\qm$ which 
allow an interpretation as ``rapidly decreasing'' functions and as 
functions with compact support on the quantum $(n,1)$-matrix ball. 
It is shown that the invariant integral is given by a generalization 
of the quantum trace. 
If an operator representation of a first order differential calculus 
over the quantum space is known, then it can be extended to the operator 
algebras of integrable functions. 
Hilbert space representations 
of $\qm$ are investigated and classified.  
Some topological aspects concerning Hilbert space representations 
are discussed. 


{\bf MSC-class:} 17B37, 47L60, 81R50 
\end{abstract}
%
%
%
%*********************************************************************
%
\section{Introduction}
%
%**********************************************************************
                                                          \label{sec-mot}
The development of quantum mechanics at the beginning of the past century 
resulted in the discovery that nuclear physics is governed by 
non-commutative quantities.
Recently, there have been made various suggestions that 
spacetime may be described by non-commutative structures at Planck scale. 
Within this approach, quantum groups might play a fundamental role. 
They can be viewed as $q$-deformations of a classical 
Lie group or Lie algebra and allow thus an interpretation as generalized 
symmetries. At the present stage, the theory is still 
in the beginning. Before constructing physical models, one has to 
establish the mathematical foundations---most important, the machineries 
of differential and integral calculus. 

In this paper, we deal with integral calculus on non-compact 
quantum spaces. The integration theory on compact quantum groups is 
well established and was mainly developed by  
S. L. Woronowicz \cite{Wo2}. He proved the existence of a 
unique normalized invariant functional (Haar functional) on 
compact quantum groups. If one turns 
to the study of non-compact quantum groups or non-compact quantum spaces, 
one faces new difficulties which do not occur in the compact case.
For instance, we do not expect that there exists a 
normalized invariant functional on the polynomial algebra of the quantum 
space. The situation is analogous to the classical theory 
of locally compact spaces, where one can only integrate functions which vanish 
sufficiently rapidly at infinity. 

Our aim is to define appropriate classes of quantized integrable 
functions for non-compact $q$-deformed manifolds. 
The ideas are similar to those in \cite{VK1}, 
%%%%(see also \cite[Section 13.4]{CP}),
where a space of finite functions was associated to the the quantum disc. 
However, our treatment will make this construction more general and 
will allow us to consider a wider class of integrable functions. Furthermore,
the invariant integral turns out to be a generalization of the well-known 
quantum trace---an observation that provides us with a rather natural proof 
of its invariance. 

Starting point of our approach will be what we call 
an operator expansion of the action. 
Suppose we are given a Hopf *-algebra $\U$ and a $\U$-module 
*-algebra $\X$ with action~$\ang$. Let $\pi:\X\rightarrow \ld$ 
be a *-rep\-re\-sen\-ta\-tion. 
(Precise definitions will be given below.)
If for any $Z\in\U$ there exists a 
finite number of 
operators $L_i,R_i\in \ld$ such that 
\begin{equation}                                     \label{oprep}
     \pi(Z\ang x)=\sum_i L_i\pi(x)R_i,\quad \ x\in\X,        
\end{equation}
then we say that we have an {\it operator expansion of the action}.
Obviously, it is sufficient to know the operators $L_i$, $R_i$ for the 
generators of $\U$. The operators $L_i$, $R_i$ are not unique as it 
can be seen by replacing $L_i$ and $R_i$ by $(-L_i)$ and $(-R_i)$.

Let us briefly outline our method of introducing an invariant integration 
theory on non-compact quantum spaces. 
Assume that $\g$ is a finite-dimensional complex semi-simple 
Lie algebra. Let $\cU_q(\g)$ denote the corresponding 
quantized  universal enveloping algebra. 
With the adjoint action $\ad_q(X)(Y):=X_{(1)}YS(X_{(2)})$, 
$\cU_q(\g)$ becomes a $\cU_q(\g)$-module (*-)algebra. 
It is a well-known fact that, for finite dimensional 
rep\-re\-sen\-ta\-tions  $\rho$ of $\cU_q(\g)$, the quantum trace formula 
$\tr_q(X):=\tr \rho(XK_{2\omega}^{-1})$, $X\in\cU_q(\g)$, 
defines an $\ad_q$-invariant linear functional on 
$\cU_q(\g)$. 
Here, the element $K_{2\omega}\in\cU_q(\g)$ is taken such that 
$K_{2\omega}^{-1}XK_{2\omega}=S^2(X)$. 

Now consider a $\cU_q(\g)$-module *-algebra $\X$ and a  
*-rep\-re\-sen\-ta\-tion $\pi:\X\rightarrow \ld$. 
In our examples, the operator expansion \rf[oprep] of the 
$\cU_q(\g)$-action on $\X$ will resemble the adjoint action. 
Furthermore, it can be extended to the *-algebra $\ld$ turning 
$\ld$ into a $\cU_q(\g)$-module *-algebra. 
The quantum trace formula suggests that we can try to define an 
invariant integral by replacing $ K_{2\omega}$ by the operator that 
realizes the operator expansion of $K_{2\omega}$ and 
taking the trace on the Hilbert space $\Hh=\bar D$. 
Since we deal with unbounded operators, this can only be done 
for an appropriate class of operators, say $\BB$. 

First of all, the generalized quantum trace should be well defined. 
Next, we wish that $\BB$ is a $\cU_q(\g)$-module *-algebra. 
This means that $\BB$ should be stable under the action defined 
by the operator expansion. If we choose $\BB$ such that 
the closures of its elements are of trace class and that
multiplying the elements of $\BB$ by any operator appearing 
in the operator expansion yields an element of $\BB$, then 
$\BB$ is certainly stable under the action of $\cU_q(\g)$ on 
$\ld$ and the generalized quantum trace is well defined on $\BB$. 
Our intention is to interpret $\BB$ as the rapidly decreasing functions 
on a $q$-deformed manifold. For this reason, we suppose additionally 
that $\BB$ is stable under multiplication by elements of $\X$. 

Clearly, the assumptions on $\BB$ are satisfied by the *-algebra 
of finite rank operators $\FF$ in $\ld$. 
The elements of $\FF$ are considered as functions with finite support 
on the $q$-deformed manifold. 

If we think of $\cU_q(\g)$ as generalized differential operators, 
then we can think of $\BB$ and $\FF$ as infinitely differentiable 
functions since both algebras are stable under the action 
of  $\cU_q(\g)$.

The algebras $\BB$ and $\FF$ were mainly introduced in order 
to develop an invariant integration theory on 
$q$-deformed manifolds. Nevertheless, our approach also 
allows to include differential calculi. 
By means of an operator rep\-re\-sen\-ta\-tion of 
a first order differential calculus 
%(abbreviated as FODC)
over $\X$, 
one can build a differential calculus over the operator 
algebras $\BB$ and $\FF$.   
In this case, we view 
the differential calculus over $\BB$ and $\FF$ as an extension of the 
differential calculus over $\X$. 

There is another notable feature of our approach. The algebras 
$\X$ (more exactly, $\pi(\X)$), $\BB$, and $\FF$ are 
subalgebras of $\ld$. In particular, they are subspaces 
of the topological space $\ldd$. Therefore we can view this 
algebras as topological spaces in a rather natural way. 
As a consequence, it makes sense to discuss topological 
concepts such as continuity, density, etc.  

In this paper, we treat the quantum $(n,1)$-matrix ball $\qm$ 
as a $\sutn$-module *-algebra \cite{VK2}.  
Since our approach to invariant integration theory  
is based on Hilbert space rep\-re\-sen\-ta\-tions, we 
shall also study *-rep\-re\-sen\-ta\-tions of $\qm$. 
When $n=1$, $\qm$ is referred to 
as quantum disc $\qd$ \cite{VK1}.  
As the algebraic relations and the *-rep\-re\-sen\-ta\-tions of $\qd$ are 
comparatively simple, it will serve as a guiding example in order to
motivate and illustrate our ideas and, therefore, we shall discuss it in a 
greater detail. 
%
%
%
%**********************************************************************
\section{Preliminaries}
%*********************************************************************
%
%++++++++++++++++++++++++++++++++++++++++++++++++++++++++++++++++++++++++
\subsection{Algebraic preliminaries}
%++++++++++++++++++++++++++++++++++++++++++++++++++++++++++++++++++++++++
                                                          \label{algprel}
Throughout this paper, $q$ stands for a real number such that $0<q<1$, 
and we abbreviate $\lambda = q-q^{-1}$.
 
Let $\cU$ be a Hopf algebra. 
The comultiplication, the counit, and the antipode of a Hopf algebra are
denoted by $\Delta$, $\varepsilon$, and $S$, respectively.
For the comultiplication $\Delta$, we employ 
the Sweedler notation: $\Delta(x)=x_{(1)}\otimes x_{(2)}$.
The main objects of our investigation are $\cU$-module algebras. An algebra 
$\X$ is called a {\it left $\cU$-module algebra} if $\X$ is a left  
$\cU$-module with action $\ang$ satisfying 
\begin{equation}                                     \label{modalg}
   f\ang (xy)=(f_{(1)}\ang x)(f_{(2)}\ang y),\quad x,y\in\X,\ f\in\cU.
\end{equation}
For an algebra $\X$ with unit $1$, we additionally require 
\begin{equation}                                      \label{modeins}
   f\ang 1=\varepsilon(f)1,\quad f\in\cU.
\end{equation}
Let $\X$ be a *-algebra and $\cU$ a Hopf *-algebra. Then $\X$ is said to 
be a {\it left $\U$-module *-al\-ge\-bra} 
if $\X$ is a left $\cU$-module algebra 
such that the following compatibility condition holds
\begin{equation}                                      \label{modstar}
   (f\ang x)^* =S(f)^*\ang x^*,\quad x\in X,\ f\in\cU.
\end{equation}

By an {\it invariant integral} we mean a linear functional $h$ on $\X$ 
such that
\begin{equation}                                       \label{invint}
   h(f\ang x)=\varepsilon(f) h(x),\quad x\in\X,\ f\in \cU.
\end{equation}
Synonymously, we refer to it as {\it $\cU$-invariant}. 

A {\it first order differential calculus} (abbreviated as FODC) over an 
algebra 
$\X$ is a pair $(\Gamma,\dd)$, where $\Gamma$ is an $\X$-bimodule and 
$\dd:\X\rightarrow \Gamma$ a linear mapping, such that 
\begin{align*}
& \dd(xy)=x\cdot\dd y+\dd x\cdot y,\quad x,y\in\X,
&  \Gamma=\Lin \{\, x{\cdot}\dd y{\cdot} z\,;\,x,y,z\in\X \,\}.
\end{align*}
$(\Gamma,\dd)$ is called a {\it first order differential *-calculus} over a 
*-algebra $\X$ if the complex vector space $\Gamma$ carries an involution *
such that 
\begin{align*}
  & (x\cdot\dd y\cdot z)^*=z^*\cdot \dd (y^*)\cdot x^*,\quad
  x,y,z\in \X.
\end{align*}

Let $(a_{ij})_{i,j=1}^n$       \label{cartan}
be the Cartan matrix of  
$sl(n+1,\C)$, that is, $a_{jj}=2$
for $j=1,\ldots,n$, $a_{j,j+1}=a_{j+1,j}=-1$ for 
$j=1,\ldots,n-1$ and $a_{ij}=0$ otherwise.
The Hopf algebra $\sltn$ is 
generated by $K_j$, $K_j^{-1}$, $E_j$, $F_j$, 
$j=1,\ldots, n$, subjected to the relations   
\begin{equation}                                 \label{sutm1}
 K_iK_j=K_jK_i,\  K_j^{-1}K_j=K_jK_j^{-1}=1,\  
 K_iE_j=q^{a_{ij}}E_jK_i,\   K_iF_j=q^{-a_{ij}}F_jK_i,
\end{equation}
\begin{equation}                                 \label{sutm2}
 E_iE_j-E_jE_i=0,\ \, i\neq j\pm1,
\quad\, 
E_j^2E_{j\pm1}-(q+q^{-1})E_jE_{j\pm1}E_j+E_{j\pm1}E_j^2=0,
\end{equation}
\begin{equation}                                \label{sutm3}
F_iF_j-F_jF_i=0,\ \,i\neq j\pm1,
\quad\,
F_j^2F_{j\pm1}-(q+q^{-1})F_jF_{j\pm1}F_j+F_{j\pm1}F_j^2=0,
\end{equation}
\begin{equation}                                \label{sutm4}
  E_iF_j-E_jF_i=0,\ \, i\neq j,  \,  \quad\,
     E_jF_j-F_jE_j=\lambda^{-1}(K_j-K_j^{-1}),\quad j=1,\ldots, n.
\end{equation}
The comultiplication $\Delta$, counit $\vare$, and antipode $S$ 
are given by
\begin{equation}            \label{delta}
\Delta (E_j)=E_j\otimes 1+K_j\otimes E_j, \quad
\Delta(F_j)=F_j\otimes K_j^{-1}+1\otimes F_j,\quad
\Delta(K_j)=K_j\otimes K_j,
\end{equation}   
\begin{equation}   
\varepsilon(K_j)=\varepsilon(K_j^{-1})=1,\quad 
\varepsilon(E_j)=\varepsilon(F_j)=0,
\end{equation}
\begin{equation}   
 S(K_j)=K_j^{-1},\quad  S(E_j)=-K_j^{-1}E_j,\quad S(F_j)=-F_jK_j.
\end{equation}
Consider the involution on $\sltn$ which is determined by 
\begin{equation}                                 \label{sutm5}    
K_i^*=K_i,\quad E_j^*=K_jF_j,\ \, F_j^*=E_jK_j^{-1},\ j\neq n,
\quad E_n^*=-K_nF_n,\ \, F_n^*=-E_nK_n^*. 
\end{equation}
The corresponding Hopf *-algebra is denoted by $\sutn$. 

If $n=1$, we write $K$, $K^{-1}$, $E$, $F$ rather than 
$K_1$, $K_1^{-1}$, $E_1$, $F_1$. These generators satisfy the 
following relations: 
\begin{equation}                                     \label{sut1}
KK^{-1}=K^{-1}K=1,\quad KEK^{-1}=q^2E,\quad KFK^{-1}=q^{-2}F, 
\end{equation}
\begin{equation}                                     \label{sut2}
        EF-FE=\lambda^{-1}(K-K^{-1}). 
\end{equation}
The involution on $\sut$ is given by  
\begin{equation}   
 K^*=K,\quad  E^*=-KF,\quad F^*=-EK^{-1}.
\end{equation}

If $n>1$, then  $K_j$, $K_j^{-1}$, $E_j$, $F_j$, 
$j=1,\ldots, n-1$ with relations \rf[sutm1]--\rf[sutm5] generate the 
Hopf *-algebra $\sutm$.
%
%
%++++++++++++++++++++++++++++++++++++++++++++++++++++++++++++++++++++++++
\subsection{Operator theoretic preliminaries}
%++++++++++++++++++++++++++++++++++++++++++++++++++++++++++++++++++++++++
                                                  \label{oppre}
We shall use the letters $\Hh$ and $\K$ to denote 
complex Hilbert spaces.
If $I$ is an at most countable index set and
$\Hh=\mathop{\oplus}_{i\in I}\Hh_i$, where $\Hh_i=\K$ for all 
$i\in I$, we denote by $\eta_i$ the vector of $\Hh$ which has the 
element $\eta\in\K$ as its $i$-th component and zero otherwise.
It is understood that $\eta_i=0$ whenever $i\notin I$. 

If $T$ is an (unbounded) operator on $\Hh$, we denote by 
$D(T)$, $\sigma(T)$, $\bar{T}$, and $T^*$ the domain, the spectrum, 
the closure, and the adjoint of $T$, respectively. 
A self-adjoint operator $A$ is called strictly positive if 
$A\ge 0$ and  ${\rm ker}\,A=\{0\}$. 
We write $\sigma (A) \sqsubseteq (a,b]$ if $\sigma(A) {\subseteq} [a,b]$ and 
$a$ is not an eigenvalue of $A$.
By definition, two self-adjoint operators strongly commute if their spectral 
projections mutually commute.

Let $D$ be a dense subspace of $\Hh$. Then the vector space 
$$
   \ld:=\{\,x\in {\rm End}(D)\,;\,D\subset D(x^*),\ x^*D\subset D\,\} 
$$
is a unital *-algebra of closeable operators with the involution 
$x\mapsto x^+:=x^*\lceil D$ and the operator product as its multiplication. 
Since it should cause no confusion, we shall continue to write 
$x^*$ in place of $x^+$. Unital *-subalgebras of $\ld$ are
called {\it O*-algebras}.

Two *-subalgebras of $\ld$ which are not O*-algebras will be of 
particular interest: The *-algebra of all finite rank operators
\begin{equation}                                  \label{FD} 
  \FF(D):=
\{\, x\in \ld\,;\, \bar{x}\ \mbox{is\ bounded},\ 
 {\rm dim}(\bar x \Hh)<\infty,\ \bar{x}\Hh\subset D,\ \bar{x}^*\Hh\subset D\,\}
\end{equation}
and, given an O*-algebra $\fA$,
\begin{equation}                                  \label{BA}
 \BB_1(\fA):=\{\,t\in \ld\,;\, \bar t\Hh\subset D,\ \bar t^*\Hh \subset D,\ 
  \ov{atb}\ \,\mbox{{\rm is\ of\ trace\ class\ for\ all}}\ a,b\in\fA\,\}. 
\end{equation}
It follows from \cite[Lemma 5.1.4]{S} that $\BB_1(\fA)$ is a *-subalgebra 
of $\ld$. Obviously, we have $\FF(D)\subset\BB_1(\fA)$ and
$1\notin \BB_1(\fA)$ if dim$(\Hh)=\infty$. 
%
%             FOLGENDES IST NOTWENDIG?
An operator 
$A\in\FF(D)$ can be written as $A=\sum^{n}_{i=1}\alpha_i e_i\otimes f_i$, 
where 
$n\in\N$, $\alpha_i\in\C$,
$f_i,e_i\in D$,
and $(e_i\otimes f_i)(x):=f_i(x)e_i$ for $x\in D$.

Assume that $\fA$ is an O*-algebra on a dense domain $D_\fA$. 
A natural choice for a topology on $D_\fA$ is the {\it graph topology} 
$t_\fA$ generated by the family of semi-norms 
\begin{equation}                                  \label{gtop} 
  \{\,||\cdot||_a\, \}_{a\in \fA},\quad ||\varphi||_a:=||a\varphi||,\quad 
  \varphi\in D_\fA.
\end{equation}
$\fA$ is called {\it closed} if the locally convex space $D_\fA$ is complete. 
The {\it closure $\bar\fA$} of $\fA$ is defined by 
\begin{equation}                                      \label{clo}
D_{\bar\fA}:=\cap_{a\in\fA}D(\bar a), \quad 
\bar\fA:=\{\,\bar a\lceil {D_{\bar{\fA}}}\,;\,a\in\fA\,\}.
\end{equation}
By \cite[Lemma 2.2.9]{S}, $D_{\bar \fA}$ is complete.

We say that $\fA$ is a {\it commutatively dominated O*-algebra on 
the Frechet domain $D_\fA$} if it satisfies the following assumptions 
(which are consequences from the definitions given in \cite{S}). There 
exist a self-adjoint operator $A$ on $\Hh$ and a sequence of Borel measurable 
real-valued functions $r_n$, $n\in\N$, such that 
$1\le r_1(t)$, $r_n(t)^2\le r_{n+1}(t)$, $r_n(A)\lceil {D_\fA}\in\fA$, 
and $D_\fA=\cap_{n\in\N} D(r_n(A))$.

Let $D_\fA^\prime$ denote the strong dual of the locally convex space $D_\fA$.
Then the conjugate space $D_\fA^+$ is the topological space $D_\fA^\prime$ 
with the addition defined as before and the multiplication replaced by 
$\alpha \cdot f:=\bar\alpha f$, $\alpha\in\C$, $f\in D_\fA^\prime$.
For $f\in D^+_\fA$ and $\varphi\in D_\fA$, we shall write 
$\ip{f}{\varphi}$ rather than $f(\varphi)$. 
The vector space of all continuous linear operators mapping 
$D_\fA$ into $D_\fA^+$ is denoted by $\ldda$. 
We assign to $\ldda$ the {\it bounded topology} $\tau_b$ generated 
by the system of semi-norms
$$
  \{\,p_M\,;\, M\subset D_\fA,\ \mbox{bounded}\,\},\quad 
  p_M(A):= {\rm sup}_{\varphi,\psi\in M}|\ip{A\varphi}{\psi}|,
  \quad A\in\fA.
$$
Notice that $\fA\subset \ldda$ for any O*-algebra $\fA$. 
Furthermore, it is known 
that $\ldop\subset \ldda$ if $D_\fA$ is a Frechet space. 

By a *-rep\-re\-sen\-ta\-tion $\pi$ of a *-algebra $\fA$ on a domain 
$D$ we mean a *-ho\-mo\-mor\-phism $\pi : \fA \rightarrow \ld$. 
For notational simplicity, we usually suppress  the rep\-re\-sen\-ta\-tion 
and write $x$ instead of $\pi(x)$ when no confusion can arise. 
If each decomposition $\pi=\pi_1\oplus\pi_2$ of $\pi$ as direct sum 
of *-rep\-re\-sen\-ta\-tions $\pi_1$ and $\pi_2$ implies that $\pi_1=0$ or
$\pi_2=0$, then $\pi$ is said to be {\it irreducible}.  

Given a *-rep\-re\-sen\-ta\-tion $\pi$, it follows from 
\cite[Proposition 8.1.12]{S} that the  mapping 
\begin{equation*}                                  %\label{crep}
     \bar\pi\,:\,\fA\rightarrow \ldop[D(\bar\pi)],
        \quad \bar\pi(a):=\ov{\pi(a)}\lceil {D(\bar\pi)}, 
\end{equation*}
defines a *-rep\-re\-sen\-ta\-tion on 
$D(\bar\pi):=\cap_{a\in\fA}D(\ov{\pi(a)})$. $\bar\pi$ is called the 
{\it closure} of $\pi$ and $\pi$ is said to be 
{\it closed} if $\bar\pi=\pi$.

If we consider *-rep\-re\-sen\-ta\-tions of *-algebras, 
we shall restrict ourself 
to rep\-re\-sen\-ta\-tions which are in a certain sense ``well behaved''. 
This means that we shall impose some regularity conditions on the 
(in general) unbounded operators under consideration. Such 
*-rep\-re\-sen\-ta\-tions will be called {\it admissible}. 
The requirements will strongly 
depend on the situation. 
Therefore there is no general definition of ``admissible''.
For further discussion on ``well behaved'' rep\-re\-sen\-ta\-tions, see 
\cite{S-hep, BIO, BIK}.

Suppose that $\X$ is a *-algebra and 
$\pi:\X\rightarrow \ld$ a *-rep\-re\-sen\-ta\-tion. 
Each symmetric operator $C\in \ld$ gives rise to a 
first order differential *-calculus $(\Gamma_{\pi,C},\dd_{\pi,C})$
over $\X$ defined by
\begin{align}
& \Gamma_{\pi,C}:=                                  \label{Gammaop}
  \Lin\{\, \pi(x)(C\pi(y)-\pi(y)C)\pi(z)\,;\,x,y,z\in\X\,\}\ 
\mbox{{\rm and}}\\
&   \dd_{\pi,C}:\X\rightarrow \Gamma_{\pi,C},\quad   \label{dop}
      \dd_{\pi,C}(x):=\im(C\pi(x)-\pi(x)C),\quad x\in\X,
\end{align}
where $\im$ denotes the imaginary unit (see \cite{S-com}).
Let $(\Gamma,\dd)$ be a first order differential *-calculus
over $\X$. Then $(\Gamma_{\pi,C},\dd_{\pi,C})$ is called 
a {\it commutator rep\-re\-sen\-ta\-tion of $(\Gamma,\dd)$}, if there 
exits a linear mapping $\rho:\Gamma\rightarrow\Gamma_{\pi,C}$ 
such that 
$\rho(x{\cdot}{\rm d} y{\cdot}z)=\pi(x)\dd_{\pi,C}(y)\pi(z)$
 and $\rho(\gamma^*)=\rho(\gamma)^*$ for all $x,y,z\in\X$, 
$\gamma\in\Gamma$.

We close this subsection by stating three auxiliary lemmas.
\begin{thl}                                     \label{L2}
Let $A$ be a  self-adjoint operator 
and let $w$ be an unitary operator on a Hilbert space $\Hh$ such that
\begin{equation}                                \label{waw}
qwA \subseteq Aw.
\end{equation}
\begin{enumerate}
\item
  The spectral projections of $A$ corresponding to $(-\infty, 0)$, $\{0\}$, and $(0, \infty)$ commute with $w$.
\item 
  Suppose additionally that  $A$ is strictly positive.
Then there exists a self-adjoint operator 
$A_0$ on a Hilbert space $\Hh_0$ 
with $\sigma (A_0) \sqsubseteq (q,1]$ such that, up to 
unitary equivalence,
$\Hh=\mathop{\oplus}^{\infty}_{n=- \infty}\Hh_n$,  
$\Hh_n=\Hh_0$, and 
$$ 
A \eta_n = q^{n} A_0 \eta_n,\quad w\eta_n = \eta_{n+1},
$$
where $\eta\in \Hh_0$ and $n\in\Z$. 
\end{enumerate}
\end{thl}
{\bf Proof.} (i):
Let $e(\mu)$ denote the spectral 
projections of $A$. Since $w$ is unitary, (\ref{waw}) implies that 
$A= q wA w^\ast$ 
and hence $e(q\mu)= we(\mu) w^\ast$. This proves (i).

\noindent
(ii): Let $\Hh_n:=e((q^{n+1},q^n])\Hh$ and 
$A_n:= A \lceil \Hh_n$, $n \in \Z$. 
Since $A$ is strictly positive,  
$\Hh=\mathop{\oplus}^{\infty}_{n=- \infty}\Hh_n$. 
Now
$e((q^{n+1},q^n])= we((q^{n},q^{n-1}]) w^\ast$
 yields $w\Hh_n =\Hh_{n+1}$.
Up to unitary equivalence, we can assume that
$\Hh_n =\Hh_0$ and $w\eta_n =\eta_{n+1}$ for $\eta \in \Hh_0$. 
Moreover,  $A\eta_n =q^n w^nAw^{n\ast}\eta_n= q^nw^nA_0\eta_0 
=q^nA_0\eta_n$. 
\hfill $\Box$
\begin{thl}                                     \label{L2a}
Let $A$ be a  self-adjoint operator and 
let $w$ be a linear isometry on a Hilbert space $\Hh$ such that
\begin{equation}                                \label{swaw}
swA \subseteq Aw
\end{equation}
for some fixed positive real number $s \neq 1$.
Suppose that $A$ has an eigenvalue $\lambda$
such that 
the eigenspace $\Hh_0:=\Ker(A-\lambda)$ coincides with 
$\Ker w^\ast$. Then the eigenspace
$\Hh_n:=\Ker(A-s^n \,\lambda)$  coincides with 
$w^n \Hh_0$ for each $n \in \N$.
\end{thl}
{\bf Proof.} 
Taking adjoints in (\ref{swaw}) gives $s^{-1} w^\ast A \subseteq A w^\ast$.
Let $n \in \N_0$, $\varphi \in \Hh_n$, and 
$\psi \in \Hh_{n+1}$. Then 
$Aw \varphi = swA \varphi = s^{n+1} \lambda w \varphi$ and
$A w^\ast \psi  = s^{-1}  w^\ast A \psi = s^n \lambda w^\ast \psi$.
Hence $w\Hh_n \subset \Hh_{n+1}$ and 
 $w^\ast\Hh_{n+1}\subset \Hh_n$.
Since $ \Hh_{n+1} \perp \Hh_0$, we have  
$w\hspace{1pt}w^\ast \psi = \psi$. This together with
$w^\ast\hspace{1pt} w =1 $ implies that $w \lceil \Hh_n$ is a bijective 
mapping from $\Hh_n$ onto $\Hh_{n+1}$ with inverse $w^\ast \lceil \Hh_{n+1}$.
\hfill $\Box$
\begin{thl}                           \label{L3}
Let $\epsilon\in\{\pm 1\}$. 
Assume that  $x$ is a closed, 
densely defined operator on a Hilbert space $\Hh$.
Then we have $\D(x x^\ast)=\D(x^\ast x)$ and the relation 
\begin{equation}                                         \label{qhyp}
x x^\ast -q^2 x^\ast x = \epsilon(1-q^2)
\end{equation} 
holds if and only if $x$ is unitarily equivalent to an orthogonal direct 
sum of operators of the following form.
\begin{description}
\item[$\epsilon=1$:]
  \begin{description}   \item[ ]
   \item[$(I)$] 
     $x\eta_n = (1-q^{2n})^{1/2} \eta_{n-1}$ on 
     the Hilbert space
     $\Hh={\oplus}^{\infty}_{n=0}\Hh_n$, $\Hh_n =\Hh_0$.
   \item[$(II)_A$] 
     $x\eta_n =(1+q^{2n} A)^{1/2} \eta_{n-1}$
     on $\Hh={\oplus}^{\infty}_{n=-\infty}\Hh_n$, $\Hh_n=\Hh_0$,
     where $A$ is a self-adjoint operator on 
     $\Hh_0$ such that $\sigma(A) \sqsubseteq (q^2,1]$.  
   \item[$(III)_u$] $x=u$, where $u$ is a unitary
     operator on $\Hh$.
  \end{description}
\item[$\epsilon=-1$:]
  \begin{description}   \item[ ]
    \item[ ] $x\eta_n = (q^{-2n}-1)^{1/2} \eta_{n+1}$ on the 
    Hilbert space  $\Hh={\oplus}^{\infty}_{n=1}\Hh_n$, $\Hh_n =\Hh_1$.
  \end{description}
\end{description}
\end{thl}
{\bf Proof.} 
Direct calculations show that  the operators described in 
Lemma \ref{L3} satisfy  (\ref{qhyp}). 
Suppose now we are given an operator 
$x$ satisfying the assumptions of the lemma.
Let $e(\mu)$ denote the spectral projections of the self-adjoint
operator $Q=\epsilon-x^\ast x$. For $\varphi \in D(Q^2)=D((x^\ast x)^2)$, 
it follows from  (\ref{qhyp}) that
\begin{align}                                         
& Qx^\ast \varphi=x^\ast (\epsilon -x x^\ast) \varphi=x^\ast 
	(\epsilon - q^2 x^\ast x - \epsilon(1-q^2)) \varphi 
	=q^2 x^\ast Q \varphi,	 \label{qxQ1}  \\
& x Q \varphi=(\epsilon -x x^\ast) x \varphi
	= (\epsilon - q^2 x^\ast x - \epsilon(1-q^2)) x \varphi
	= q^2 Q x   \varphi. 	 \label{qxQ2}  
\end{align}
The cases $\epsilon=1$ and $\epsilon=-1$ will be analyzed separately.

\noindent
$\epsilon=1$: Let $x^\ast =u a$ be the polar decomposition
of $x^\ast$. Note that 
\begin{equation}                                         \label{b1}
a^2=x x^\ast = 1-q^2 + q^2 x^\ast x = 1-q^2 Q \geq 1-q^2 ,
\end{equation} 
which implies, in particular, that $\Ker a= \Ker u =0$, so $u$ is an
isometry. Inserting $\varphi=a^{-1} \psi$ in
(\ref{qxQ1}), where $\psi \in D(Q^{3/2})$, one obtains 
$Q u \psi = q^2 u a Q a^{-1} \psi = q^2 u Q \psi$.
Since $D(Q^{3/2})$ is a core for $Q$, it follows that $q^2 u Q \subseteq Q u$.
By taking adjoints, one also gets $u^\ast Q \subseteq q^2 Q u^\ast$. 
Furthermore, 
$\varphi \in \Ker x= \Ker x^\ast x = \Ker u^\ast$ if and only if 
$(Q-1) \varphi =0$. If $\Ker u^\ast \neq \{ 0 \}$, 
Lemma \ref{L2a} implies that 
$\K:=\mathop{\oplus}_{n=0}^\infty \, \Hh_n$, where 
$\Hh_n= \Ker (Q\hspace{1.5pt}{-}\hspace{1.5pt}q^{2n})$,
is a reducing subspace for $u$ and $Q$. Moreover,
$x\lceil \K=(1\hspace{1.5pt}{-}\hspace{1.5pt}q^2 Q)^{1/2} u^\ast \lceil \K$ 
is unitarily equivalent to an operator of the form $(I)$.

It suffices now to prove the assertion under the additional assumption that
$\Ker u^\ast = \{0\}$. By Lemma \ref{L2}(i), 
we can treat the cases where 
$Q$ is strictly positive, zero, or strictly negative separately.

If $Q$ were strictly positive , then  it would be  unbounded by 
Lemma \ref{L2}(ii), which contradicts (\ref{b1}).  
Hence we can discard this case.
If $Q=0$, then $x=u^\ast$ is unitarily equivalent to an operator of
the form $(III)_u$.
When  $Q$ is strictly negative, Lemma \ref{L2}(ii)  applied to the relation
$q^2 u (-Q) \subseteq (-Q) u$
shows  that $x=(1-q^2  Q)^{1/2} u^\ast$  is unitarily equivalent to an 
operator of the form $(II)_A$.

\noindent
$\epsilon=-1$: In this case,
we use  the polar decomposition $x=v b$ of $x$. From 
\begin{equation}                                         \label{b-1}
b^2=x^\ast x = -1-Q= q^{-2}(x x^\ast +1 -q^2)  \geq q^{-2} -1,
\end{equation} 
it follows that $\Ker b= \Ker v =\{0\}$ so that $v$ is an
isometry. Using (\ref{qxQ2}) and arguing as above, 
one obtains $q^{-2} v Q \subseteq Q v$ and 
$q^2 v^\ast Q \subseteq Q v^\ast$. 
Note that, in the present case, $Q \leq -q^{-2}$  
by (\ref{b-1}). Therefore 
$\Ker v^\ast \neq \{0\}$ since otherwise Lemma \ref{L2} would imply that 
$0$ belongs to the spectrum of $Q$.  
Now $\varphi \in \Ker v^\ast=\Ker x^\ast =\Ker x x^\ast $ if and only if
$Q \varphi=(-1-x^\ast x) \varphi = (-1-q^{-2}(1-q^2)) \varphi 
= -q^{-2} \varphi$.
From Lemma \ref{L2a}, it follows that 
$\K:=\mathop{\oplus}_{n=1}^\infty  \Hh_n$, where $\Hh_n= \Ker (Q+q^{-2n}$),
is a reducing subspace for $v$ and $Q$. In particular, 
$x\lceil \K =v (-1-Q)^{1/2} \lceil \K$  is unitarily equivalent to 
an operator of the form stated in the lemma.  
Finally, we conclude that $\Hh=\K$ since the restriction of $v^\ast$ to 
a nonzero 
orthogonal complement of $\K$ would be injective, 
which is impossible as noted before.  
\hfill $\Box$
\mn

\noindent
{\bf Remark.}
For $\epsilon=1$, 
a characterization of irreducible representations of (\ref{qhyp})
can be found in \cite{PW} as a special case of the results therein.  
For $\epsilon=-1$, the irreducible representations of (\ref{qhyp})
were obtained in  \cite{CGP} by assuming in the proof that $x^\ast x$ has eigenvectors.
%
%
%++++++++++++++++++++++++++++++++++++++++++++++++++++++++++++++++++++++++
\section{Quantum disc}
%++++++++++++++++++++++++++++++++++++++++++++++++++++++++++++++++++++++++
                                                             \label{QD}
%
%------------------------------------------------------------------------
\subsection{Invariant integration on the quantum disc}
%------------------------------------------------------------------------
                                                           \label{sec-qd}
The quantum disc $\qd$ is defined as the *-algebra generated 
by $z$ and $z^*$ with relation
\begin{equation}                                   \label{qd-rel}
               z^*z-q^2zz^*=1-q^2.
\end{equation}
By \rf[qd-rel], it is obvious that $\qd=\Lin\{ z^nz^{*m}\,;\,n,m\in\N\}$.
Set 
\begin{equation}                                   \label{y}
                          y:=1-zz^*.
\end{equation}
Then $y=y^*$ and 
\begin{equation}                                   \label{yz}
        yz=q^2zy,\quad yz^*=q^{-2}z^*y.
\end{equation}
From $zz^*=1-y$, $z^*z=1-q^2y$, and \rf[yz], we deduce
\begin{equation}                                   \label{zz*}
        z^nz^{*n}=(y;q^{-2})_n,\quad z^{*n}z^n=(q^2y;q^2)_n,
\end{equation}
where $(t;q)_0:=1$ and  $(t;q)_n:=\prod_{k=0}^{n-1}(1-q^kt)$, $n\in\N$. 
In particular, each element $f\in \qd$ can be written as 
\begin{equation}                                   \label{f=}
   f=\sum_{n=0}^N z^np_n(y)+\sum_{n=1}^M p_{-n}(y)z^{*n},
        \quad N,M\in\N,
\end{equation}
with polynomials $p_n$ in $y$.

The left action $\ang$ which turns $\qd$ into a $\sut$-module *-algebra can 
be found in \cite{VK1, VK3} or \cite{KL}.  On generators, it takes the 
form 
\begin{equation}                                     \label{qd-act}
 K^{\pm1}\ang z=q^{\pm 2}z,\quad E\ang z=-q^{1/2}z^2,\quad 
 F\ang z=q^{1/2},
\end{equation}
\begin{equation}                                     \label{qd-act*}
 K^{\pm}\ang z^*=q^{\mp2}z^*,\quad E\ang z^*=q^{-3/2},\quad 
 F\ang z^*=-q^{5/2}z^{*2}.
\end{equation}

Remind our notational conventions regarding 
rep\-re\-sen\-ta\-tions. For instance, if $\pi:\qd\rightarrow \ld$ is a 
rep\-re\-sen\-ta\-tion, we write $f$ instead of $\pi(f)$ and 
$X\ang f$ in instead of $\pi(X\ang f)$, where $f\in\qd$, $X\in\sut$. 
The key observation of this subsection is the following simple 
operator expansion.
\begin{thl}                                           \label{discact}
Let $\pi:\qd\rightarrow \ld$ be a *-rep\-re\-sen\-ta\-tion of $\qd$ such 
that $y^{-1}$ belongs to $\ld$. 
Set $A:=q^{-1/2}\lambda^{-1}z$ and $B:=-y^{-1}A^*$. Then the formulas 
\begin{align}
 & K\ang f =yfy^{-1},\quad K^{-1}\ang f=y^{-1}fy,\label{opact1}\\
 & E\ang f=Af-yfy^{-1}A,                          \label{opact2}\\
 & F\ang f= Bfy-q^2fyB                 \label{opact3}     
\end{align}       
define an operator expansion of the action $\ang$, where $f\in\qd$. 
The same formulas applied to $f\in \ld$ turn the 
$O^*$-algebra $\ld$ into a $\sut$-module *-algebra.
\end{thl}

{\bf Proof.} We take Equations \rf[opact1]--\rf[opact3] as definition and 
show that the action $\ang$ defined in this way turns $\ld$ into a 
$\sut$-module *-algebra. 
To verify that $\ang$ is well defined, we use the commutation relations 
\begin{equation}                                     \label{yT}
   yA=q^2Ay,\quad yB=q^{-2}By,\quad AB-BA=-\lambda^{-1}y^{-1}
\end{equation}
which are easily obtained by applying \rf[qd-rel] and \rf[yz]. 
Let $f\in \ld$. It follows that 
$$                               
KE\ang f=y(Af-yfy^{-1}A)y^{-1} =q^2(Ayfy^{-1}-y^2fy^{-2}A)=q^2EK\ang f
$$
and 
\begin{eqnarray*}                                   %  \label{}
(EF-FE)\ang f &=&ABfy+yfBA-BAfy-yfAB\\
              &=&(AB-BA)fy-yf(AB-BA) \\
                  &=&\lambda^{-1} (yfy^{-1}-y^{-1}fy)
             \ =\ \lambda^{-1}(K-K^{-1})\ang f.
\end{eqnarray*}
The other relations of \rf[sut1] are treated in the same way, so we conclude 
that the action is well defined. 

We continue by verifying \rf[modalg]--\rf[modstar]. Since the action is 
associative, it is sufficient to prove \rf[modalg]--\rf[modstar] for
the generators $E$, $F$, $K$, and $K^{-1}$. 
Let $f,g\in\ld$. Then
$$
K^{\pm1}\ang (fg)=y^{\pm1}fgy^{\mp1}=y^{\pm1}fy^{\mp1}y^{\pm1}gy^{\mp1}
                 =(K^{\pm1}\ang f)(K^{\pm1}\ang g).
$$
Furthermore,
\begin{eqnarray*}
(E\ang f)g+(K\ang f)(E\ang g) 
           &=& (Af-yfy^{-1}A)g +yfy^{-1}(Ag-ygy^{-1}A)\\
           &=&  Afg-yfgy^{-1}A \ =\  E\ang (fg)
\end{eqnarray*}
and, analogously, $(F\ang f)(K^{-1}\ang g)+f(F\ang g)=F\ang (fg)$. 

Clearly, $K^{\pm1} \ang 1=\vare(K^{\pm1})1$,
$E\ang 1=A-yy^{-1}A=0=\vare(E)1$, and, similarly, 
$F\ang 1=0=\vare(F)1$.  
Equation \rf[modstar] reads for $K^{\pm1}$ 
$$
(K^{\pm1}\ang f)^*=(y^{\pm1}fy^{\mp1})^*=
y^{\mp 1}f^*y^{\pm1}=K^{\mp1}\ang f^*=
S(K^{\pm1})^*\ang f^*,
$$
and, for $E$, we have 
\begin{equation}                               \label{Efstar}
 (E\ang f)^*= f^*A^*-A^*y^{-1}f^*y=-f^*yB+q^{-2}Bf^*y =
  q^{-2}F\ang f^* =  S(E)^*\ang f^*
\end{equation}
since $S(E)^*=KFK^{-1}=q^{-2}F$. 
Replacing in \rf[Efstar] $f$ by $q^2f^*$ and applying the involution 
gives 
$(F\ang f)^* =q^2E\ang f^*=S(F)^*\ang f^*$, 
where we used $S(F)^*=q^2E$. 
Summarizing, we have shown that the action $\ang$ defined by 
\rf[opact1]--\rf[opact3] equips $\ld$ with the structure of a 
$\sut$-module *-algebra. 

It remains to prove that \rf[opact1]--\rf[opact3] define an 
operator expansion of the action $\ang$ given 
by \rf[qd-act] and \rf[qd-act*]. 
Since $\pi(\qd)$ is a 
*-subalgebra of the $\sut$-module *-algebra $\ld$, 
it is 
sufficient to verify \rf[opact1]--\rf[opact3] for the 
generators of $\sut$ and $\qd$ (see Equation \rf[modalg]).   
From the definition of $A$ and $y$, it follows
by using \rf[yT] and \rf[qd-rel] that 
\begin{equation}                                     \label{Kz}
K^{\pm1}\ang z=y^{\pm1}zy^{\mp1}=q^{\pm2}z,\quad 
K^{\pm1}\ang z^*=y^{\pm1}z^*y^{\mp1}=q^{\mp2}z^*,
\end{equation}
\begin{equation}                                      \label{Ez}
E\ang z=Az-yzy^{-1}A=q^{-1/2}\lambda^{-1}(z^2-q^2z^2)=-q^{1/2}z^2,
\end{equation}
\begin{equation}                                       \label{Ez*}
E\ang z^*=Az^*-yz^*y^{-1}A=q^{-5/2}\lambda^{-1}(q^2zz^*-z^*z)=q^{-3/2}
\end{equation}
and, similarly, $F\ang z=q^{1/2}$, $F\ang z^*=-q^{5/2}z^{*2}$. 
This completes the proof.                               \hfill $\Box$
\mn

Recall that the left adjoint action 
$\ad_{L}(a)(b):=a_{(1)}bS(a_{(2)})$, $a,b\in\sut$, turns 
$\sut$ into a $\sut$-module *-algebra. 
For the generators  $E$, $F$, and $K$, we obtain 
$\ad_{L}(E)(b)=Eb-KbK^{-1}E$,
$\ad_{L}(F)(b)= FbK-q^2bKF$, 
and $\ad_{L}(K)(b)=KbK^{-1}$. There is an obvious 
formal coincidence of this formulas with \rf[opact1]--\rf[opact3] but 
$A$, $B$, and $y$ do not satisfy the relations of $E$, $F$, and $K$ 
because the last equation of \rf[yT] differs from \rf[sut2]. 

We mentioned that for a finite dimensional 
rep\-re\-sen\-ta\-tion $\rho$ of $\sut$ the quantum trace 
$$
           \tr[q]a:= \tr\rho(aK^{-1})
$$ 
defines an invariant integral on $\sut$ 
(see \cite[Proposition 7.1.14]{KS}). 
The proof does not involve the whole set of relations of 
$\sut$ but the trace property and the relation 
$K^{-1}fK=S^2(f)$ for all $f\in \sut$. The last relation reads on 
generators as $K^{-1}KK=K$, $K^{-1}EK=q^{-2}E$, $K^{-1}FK=q^2F$ 
and these equations 
are also satisfied if we replace $K$ by $y$, $E$ by $A$, and $F$ 
by $B$. 

The main result of this section, achieved in Proposition \ref{haar}
below, is a generalization of the quantum trace formula 
to the operator algebras $\BBA$ and $\FFD$ from Subsection \ref{oppre} 
by using the above observations. 
Notice that we cannot have a normalized invariant integral on 
$\qd$; if there were an invariant integral $h$ on $\qd$ 
satisfying $h(1)=1$, then we would obtain
\begin{equation}                                 \label{h=0}
  1=h(1)=q^{-1/2}h(F\ang z)=q^{-1/2}\vare(F)h(z), 
\end{equation} 
a contradiction since $\vare(F)=0$.
\begin{thp}                                       \label{haar}
Suppose that $\pi:\qd\rightarrow \ld$ is a *-rep\-re\-sen\-ta\-tion of $\qd$ 
such that $y^{-1}\in\ld$. Let $\fA$ be the O*-algebra generated 
by the operators $z$, $z^*$, and $y^{-1}$.
Then the *-algebras $\FFD$ and $\BBA$ defined in 
\rf[FD] and \rf[BA], respectively, are $\sut$-module *-algebras,
where the action is given by \rf[opact1]--\rf[opact3]. 
The linear functional
\begin{equation}                                  \label{h}
    h(g):=c\, \tr \ov{gy^{-1}},\quad c\in\R,        
\end{equation}
defines an invariant integral on both $\FFD$ and $\BBA$. 
\end{thp}
{\bf Proof.}
Obviously, by the definition of $\FFD$ and $\BBA$, we have 
$afb\in\FFD$ and $agb\in\BBA$ for all $f\in\FFD$, $g\in\BBA$, 
$a,b\in\fA$, so both algebras are stable under the action of 
$\sut$. 
By Lemma \rf[discact], this action turns $\FFD$ and $\BBA$ into
$\sut$-module *-algebras. 

The proof of the invariance of $h$ uses the trace property 
$\tr \ov{agb}=\tr \ov{gba}=\tr \ov{bag}$ which holds for all 
$g\in\BBA$ and all $a,b\in\fA$ (see \cite{S}). Since the action 
is associative and $\vare$ a homomorphism, we only have to prove the 
invariance of $h$ for generators. Let $g\in\BBA$. Clearly, 
$$
h(K^{\pm1}\ang g)=\tr\ov{y^{\pm1}gy^{\mp1}y^{-1}}=
 \tr\ov{gy^{-1}}=\vare(K^{\pm1})h(g),
$$
$$
h(E\ang g)= \tr(\ov{Agy^{-1}}-\ov{ygy^{-1}Ay^{-1}})=
\tr\ov{Agy^{-1}}-\tr\ov{Agy^{-1}}=0=\vare(E)h(g).
$$
Using the second relation of 
\rf[yT], we compute 
$$
h(F\ang g)= \tr(\ov{Bg}-q^2\ov{gyBy^{-1}})=
\tr\ov{gB}-\tr\ov{gB}=0=\vare(B)h(g).
$$
Hence $h$ defines an 
invariant integral on $\BBA$. It is obvious that the restriction 
of $h$ to $\FFD$ gives an invariant integral on $\FFD$. 
 \hfill $\Box$
\mn

Commonly, the algebra $\qd$ is considered as the polynomial functions 
on the quantum disc. Observe that $agb\in\BBA$ for all
$g\in\BBA$ and all polynomial functions $a,b\in\qd$. 
Notice, furthermore, that the action of $E$ and $F$ satisfies 
a ``twisted'' Leibniz rule. 
If we think of $\sut$ as an algebra of ``generalized differential 
operators'', then we can think of $\BBA$ as the algebra of infinitely 
differentiable functions which vanish sufficiently rapidly at 
``infinity'' and of $\FFD$ as the infinitely differentiable 
functions with compact support. 
%
%
%------------------------------------------------------------------------
\subsection{Topological aspects of *-representations}
%------------------------------------------------------------------------
%
\label{topasp}
This subsection is concerned with some topological aspects of the 
rep\-re\-sen\-ta\-tions of $\qd$. 
The rep\-re\-sen\-ta\-tions of $\qd$ are given 
by Lemma \ref{L3}. Here we restate Lemma \ref{L3} by considering 
only irreducible *-rep\-re\-sen\-ta\-tions and specifying the domain on which 
the operators act. As we require that $y^{-1}$ exists, we exclude the 
case $(III)_u$ in which $y=0$. 
Let $\{\eta_j\}_{j\in J}$ denote the canonical basis in the Hilbert 
space $\Hh=l_2(J)$, where $J=\N_0$ or $J=\Z$. 
\begin{description}
\item[$(I)$] 
The operators $z$, $z^*$, and $y$ act on 
$D:=\Lin \{\,\eta_n\,;\, n\in \N_0\,\}$ by 
$$
z\eta_n=\lambda_{n+1} \eta_{n+1},\quad
z^*\eta_n=\lambda_n\eta_{n-1},\quad y\eta_n=q^{2n}\eta_n. 
$$
\item[$(II)_\alpha$]
Let $\alpha\in [0,1)$. 
The actions of $z$, $z^*$, and $y$ on 
$D:=\Lin \{\,\eta_n\,;\, n\in \Z\,\}$
are given by 
$$
z\eta_n=\lambda_{\alpha,n+1}\eta_{n+1},
\quad z^*\eta_n=\lambda_{\alpha,n}\eta_{n-1},\quad 
y\eta_n=-q^{2(\alpha+n)}\eta_n,
$$
\end{description}
where $\lambda_n=(1-q^{2n})^{1/2}$ and  
$\lambda_{\alpha,n}=(1+q^{2(\alpha+n)})^{1/2}$.
Obviously, $y^{-1}\in \ld$ in both cases. 

Let $\fA$ be the O*-algebra defined in Proposition \ref{haar}. 
If we equip $D$ with the graph topology $t_\fA$, 
$D$ is not complete. The situation becomes better
if we pass to the closure of $\fA$.
Since this can always be done, there is no loss of generality
in assuming $\fA$ to be closed, that is,  
$D_{\fA}:=\cap_{a\in\fA}D(\bar a)$ (see Equation \rf[clo]). 
Some topological facts concerning $\fA$ and $\ldop$
are collected in the following lemma and the next proposition. 
\begin{thl}                    \label{l-top}
Suppose that we are given an irreducible *-rep\-re\-sen\-ta\-tion of type 
$(I)$ or $(II)_\alpha$ and the O*-algebra $\fA$ from 
Proposition \ref{haar} is closed.  
\begin{enumerate}
\item
  $\fA$ is a commutatively dominated O*-algebra on a 
  Frechet domain.
\item 
   $D_\fA$ is nuclear, in particular, $D_\fA$ is 
   a Frechet--Montel space.
\end{enumerate}
\end{thl}
{\bf Proof.} (i): The operator $y$ is essentially self-adjoint on
$D_\fA$ and so is 
\begin{equation}                                      \label{opT}
       T:= 1+y^2+y^{-2}.
\end{equation}
Let $\varphi \in D_\fA$.
A standard argument shows that, for each polynomial $p(y,y^{-1})$,  
there exist  $k\in\N$ such that 
$||p(y,y^{-1})\varphi||\le ||T^k\varphi||$. By using \rf[zz*], we
get the estimates 
$$
||z^np(y,y^{-1})\varphi|| \le (|| \bar p(y,y^{-1})(q^2y;q^{2})_n
\, p(y,y^{-1})\varphi||\,||\varphi||)^{1/2}\le ||T^l\varphi||,
$$
 $$
||z^{*n}p(y,y^{-1})\varphi|| \le (|| \bar p(y,y^{-1})(y;q^{-2})_n
\, p(y,y^{-1})\varphi||\,||\varphi||)^{1/2}\le ||T^{l^\prime}\varphi||
$$
for some $l,l^\prime\in\N$. 
Since $T\ge 2$ and  $T^k\leq T^m$ for $k\leq m$,
we can find for each finite sequence $k_1, \ldots,k_N\in\N$ 
a $k_0\in\N$ such that
$\sum_{j=1}^N||T^{k_j}\varphi||\le ||T^{k_0}\varphi||$.  
By \rf[f=], \rf[yz], and the definition of $\fA$,  
it follows that each $f\in\fA$ can be written as  
$f=\sum_{n=0}^N z^np_n(y,y^{-1})+\sum_{n=1}^M z^{*n}p_{-n}(y,y^{-1})$. 
From the foregoing,  we conclude  that there exist $m\in\N$ 
such that $||f\varphi||\le ||T^m\varphi||$, 
consequently $||\cdot||_f\le||\cdot||_{T^m}$.
This shows that the family 
$\{ ||\cdot||_{T^{2^k}}\}_{k\in\N}$ generates the graph topology and 
$D_\fA=\cap_{k\in\N}D(\bar{T}^{2^k})$, which proves (i).
\sn\\
(ii): By (i), the graph topology is metrizable. It follows from 
\cite[Proposition 2.2.9 and Corollary 2.3.2.(ii)]{S} that 
$D_\fA$ is a reflexive Frechet space, in particular, 
$D_\fA$ is barreled. 
To see that $D_\fA$ is nuclear, consider 
$E_n:=(\ov{D_\fA},||\cdot||_{T^n})$, where the closure of $D_\fA$ 
is taken in the norm $||\cdot||_{T^n}$, and the embeddings 
$\iota_{n+1}:E_{n+1}\rightarrow E_n$, where $\iota_{n+1}$ 
denotes the identity on $E_{n+1}$, $n\in\N$. 
It is easy to see that the operator $\bar{T}^{-1}:\Hh\rightarrow\Hh$ 
is a Hilbert--Schmidt operator and that the canonical basis 
$\{ e_j\}_{j\in J}$, where $J=\N_0$ in case $(I)$ and 
$J=\Z$ in case $(II)$, is a complete set of eigenvectors. 
The set $\{ f^n_j\}_{j\in J}$, $f_j^n=||T^ne_j||^{-1}e_j$ 
constitutes an orthonormal basis in $E_n$, and we have 
\begin{align*}
 \sum_{j\in J}||\iota_{n+1}(f^{n+1}_j)||_{T^n}^2
  &=  \sum_{j\in J}||T^nf^{n+1}_j||^2
 =\sum_{j\in J}||T^{n}(||T^{n+1}e_j||^{-1}e_j)||^2 \\
  &=  \sum_{j\in J}||T^{-1}e_j||^2<\infty
\end{align*} 
which shows that $\iota_{n+1}$ is a Hilbert--Schmidt operator. 
From this, we conclude that $D_\fA$ is a nuclear space since 
the family $\{||\cdot||_{T^n}\}_{n\in\N}$ of Hilbert semi-norms
generates the topology on $D_\fA$. As each nuclear space is a 
Schwartz space and as each barreled Schwartz space is a Montel space, 
$D_\fA$ is a Montel space.
                                          \hfill $\Box$
\begin{thp}                    \label{top}
Suppose we are given an irreducible *-rep\-re\-sen\-ta\-tion of type 
$(I)$ or $(II)_\alpha$. Assume that the O*-algebra $\fA$ from 
Proposition \ref{haar} is closed.  
\begin{enumerate}
\item
   $\FFDA$ is dense in $\ldda$ with respect to the bounded 
   topology $\tau_b$.
\item
   The $\sut$-action on $\ldop$ is continuous with respect to $\tau_b$.  
\end{enumerate}
\end{thp}
{\bf Proof.}
(i) follows immediately from Lemma \ref{l-top}(ii) 
and \cite[Theorem 3.4.5]{S}. \\
(ii): Let $x\in\ldop$ and $a,b\in\fA$. 
According to \cite[Proposition 3.3.4(ii)]{S}, the mul\-ti\-pli\-ca\-tion 
$x\mapsto axb$ is continuous. By Lemma \ref{discact}, 
the action of $\sut$ is given by a finite linear combination of such 
expressions, hence it is continuous. 
                                                        \hfill $\Box$
\mn

The algebra $\FFD$ is the linear span of operators 
$\eta_m\otimes\eta_n$, where $n,m\in\N_0$ for the type $(I)$ 
rep\-re\-sen\-ta\-tion 
and $n,m\in\Z$ for type $(II)$ rep\-re\-sen\-ta\-tions. 
Since $D\subset D_\fA$, we can consider $\FFD$ as a subalgebra of 
$\FFDA$ and, moreover, as a $\sut$-module *-algebra. 
The interest in $\FFD$ stems from the fact that the operators 
$\eta_n\otimes\eta_m$ are more suitable for calculations.  With a little 
extra effort, we can deduce from Proposition \ref{top} that the linear span 
of this operators is dense in $\ldda$. 

\begin{thc}                                           \label{dense}
 $\FFD$ is dense in $\ldda$ with respect to the bounded topology $\tau_b$. 
\end{thc}   
{\bf Proof.} In view of Proposition \ref{top}(i), it is sufficient to 
show that $\FFDA$ lies in the closure of\ $\FFD$. 
With $T$ defined in \rf[opT], 
consider the set of Borel measurable functions 
$$
\cR:=\{\,r:\sigma(\bar T)\rightarrow [0,\infty)\,;\, 
\sup_{t\in\sigma(\bar T)} r(t)t^{2^n}<\infty\,\}. 
$$
It follows from Lemma \ref{l-top}(i) and  
\cite[Proposition 3.4]{K} that the family of semi-norms 
$$
\{ ||\cdot||_r\}_{r\in\cR},\quad ||a||_r:= ||\ra a\,\ra||,
  \quad a\in\ldda, 
$$
(the norm $||\cdot||$ being the operator norm in $\mathcal{L}(\Hh)$)
generates the topology $\tau_b$. 

Let $\varphi,\psi\in D_\fA$. 
Notice that 
$||\ra (\varphi\otimes\psi) \ra||\leq ||\ra||^2||\varphi||\,||\psi||$. 
With $\alpha_n,\beta_n\in\C$, write 
$\varphi=\sum_{n\in J}\alpha_n\eta_n$, 
$\psi=\sum_{n\in J}\beta_n\eta_n$,  
where $J=\N_0$ or $J=\Z$ 
according to the type of rep\-re\-sen\-ta\-tion considered. 
For $k\in\N$, set 
$\varphi_k:=\sum_{|n|\leq k}\alpha_n\eta_n$ and 
$\psi_k:=\sum_{|n|\leq k}\beta_n\eta_n$. Clearly, 
$\varphi_k,\psi_k\in\FFD$. Now
\begin{align*}
  ||\varphi\otimes\psi-\varphi_k\otimes\psi_k||_r
  &= ||\ra (\varphi\otimes\psi-\varphi_k\otimes\psi_k)\ra||\\
  &\leq ||\ra||^2||\varphi-\varphi_k||\,||\psi||
      +||\ra||^2||\varphi_k||\,||\psi-\psi_k||\rightarrow 0
\end{align*}
as $k\rightarrow \infty$ 
for all $r\in\cR$, hence $\varphi\otimes\psi$ lies in the closure 
of $\FFD$. Since $\FFDA$ is the linear span of operators 
$\varphi\otimes\psi$, the assertion follows. 
                                                       \hfill $\Box$
\mn

Proposition \ref{top} and Corollary \ref{dense} show how $\FFD$ and 
$\FFDA$ are related to the image of  $\qd$ 
in $\ldop$: By density and continuity, 
$\FFD$ and $\FFDA$ carry the whole information about the action of 
$\sut$ on $\ldop\subset\ldda$ and, in particular, on $\qd\subset\ldda$.
 
It would be desirable to have also the converse statement,
that is, to obtain the action on $\FFD$ (or $\FFDA$) by taking the 
closure of $\qd$ in $\ldda$. Unfortunately, this is not possible. 
From \cite[Theorem 4.5.4]{S}, it follows that $\tau_b$ coincides 
with the finest locally convex topology on $\fA$. 
Since $\fA$ is closed with respect to the finest locally convex
topology, it is closed with respect to $\tau_b$. 

For $\FFD$ to be in the closure of $\qd$, we can consider 
a different locally convex topology on $D$. Let $D^\prime$ 
be the vector space of all formal series 
$\sum_{j\in J}\alpha_j\eta_j$, where $J=\N$ or $J=\Z$. 
There exists a dual pairing $\ip{\cdot}{\cdot}$ 
of  $D^\prime$ and $D$ given by 
$$
\ip{\sum_{j\in J}\alpha_j\eta_j}{\sum_{|n|\leq n_0}\beta_n\eta_n}
=\sum_{|n|\leq n_0}\alpha_n\beta_n. 
$$
We equip $D$ and $D^\prime$ 
with the weak topologies arising from this dual pairing. 
To $\lddp$, the vector space of all continuous linear mappings 
from $D$ into $D^\prime$, we assign the operator weak topology 
$\tau_{ow}$, that is, the topology generated by the family of 
semi-norms 
$$
\{p_{\varphi,\psi}\}_{\varphi,\psi\in D},\quad
 p_{\varphi,\psi}(a):=|\ip{a\varphi}{\psi}|,\quad 
a\in\lddp.
$$
Then $\qd$ is dense in $\lddp$ with respect to $\tau_{ow}$ and the 
action of $\sut$ on $\ld$ defined by \rf[opact1]--\rf[opact3] 
is continuous. This is essentially the method for constructing 
the space $D(U_q)^\prime (=\lddp)$ of distributions 
on the quantum disc as performed in \cite{VK1}.
The topological space $D(U_q)$ of finite functions on the quantum disc 
defined in \cite{VK1} is homeomorphic to $\FFD$ with the 
operator weak topology $\tau_{ow}$. 

We now give another description of $\FFD$. 
\begin{thl}                                   \label{FS}
Let $\FS$ be the set of (Borel measurable) functions on 
$\sigma(\bar y)$ with finite support, that is, 
$$
  \FS=\{\, \psi : \sigma(\bar y)\rightarrow \C\,;\, 
  \#\{t\in\sigma(\bar y);\psi(t)\neq0\}<\infty\,\}. 
$$
Each $f\in \FFD$ can be written as
$$
  f=\sum_{n=0}^N z^n\psi_n(\bar y)+
     \sum_{n=1}^M\psi_{-n}(\bar y)z^{*n},
  \quad N,M\in\N,
$$
where $\psi_k\in\FS$, $k=-M,\ldots,N$. 

Conversely, if 
$\psi_k\in\FS$, then 
$\sum_{n=0}^N z^n\psi_n(\bar y)+
          \sum_{n=1}^M\psi_{-n}(\bar y)z^{*n}\in\FFD$.
\end{thl}
{\bf Proof.} To see this, consider the functions 
\[ \delta_{k}(t)\ :=\ \left\{ \begin{array}{r@{\quad :\quad}l}
                  1      & \mathrm{for}\ \ t=q^{2k}   \\
                  0      & \mathrm{for}\ \ t\neq q^{2k} 
                \end{array} \right.    \]
if we are given a type $(I)$ rep\-re\-sen\-ta\-tion, and  
\[ \delta_{k}(t)\ :=\ \left\{ \begin{array}{r@{\quad :\quad}l}
         1      & \mathrm{for}\ \ t=-q^{2\alpha+2k}   \\
         0      & \mathrm{for}\ \ t\neq -q^{2\alpha+2k} 
                \end{array} \right.    \]
if we are given a rep\-re\-sen\-ta\-tion of type $(II)_\alpha$. 
Notice that $\dy$ is the projection on $\Hh$ with range
$\C\eta_k$, that is, 
$\dy=\eta_k\otimes\eta_k$. 

Each $\psi_n\in\FS$ can be written as a finite sum
$\sum_k \psi_{n,k}\delta_k(t)$, where 
$\psi_{n,k}=\psi_n(q^{2k})$ for the type $(I)$ rep\-re\-sen\-ta\-tion 
and $\psi_{n,k}=\psi_n(-q^{2\alpha+2k})$ for type $(II)_\alpha$ 
rep\-re\-sen\-ta\-tions. Furthermore, we have 
$$
z^n\dy=z^n(\eta_k\otimes\eta_k)=(z^n\eta_k)\otimes\eta_k\in\FFD,
$$
$$
\dy z^{*n}=(\eta_k\otimes\eta_k)z^{*n}
          =\eta_k\otimes(z^n\eta_k)\in\FFD,
$$
hence $\sum_{n=0}^N z^n\psi_n(\bar y)+
\sum_{n=1}^M\psi_{-n}(\bar y)z^{*n}\in\FFD$
whenever $\psi_n\in\FS$, $n=-M,\ldots,N$. 

On the other hand, for $k\leq n$, we can write 
$$
\eta_n\otimes\eta_k=\gamma_{n,k}(z^{n-k}\eta_k)\otimes\eta_k
   =\gamma_{n,k}z^{n-k}\dy,
$$
$$
\eta_k\otimes\eta_n=\gamma_{n,k}\eta_k\otimes(z^{n-k}\eta_k)
   =\gamma_{n,k}\dy z^{*n-k},
$$
where $\gamma_{n,k}=(q^{2(k+1)};q^2)_{n-k}^{-1/2}$ and
$\gamma_{n,k}=(-q^{2(\alpha+k+1)};q^2)_{n-k}^{-1/2}$
for the rep\-re\-sen\-ta\-tions of type $(I)$ and type $(II)_\alpha$,
respectively. Hence any linear combination of 
$\eta_m\otimes\eta_l$ is equivalent to a linear combination 
of $z^n\dy$ and $\dy z^{*n}$. 
Summing over equal powers of $z$ and $z^*$ yields coefficients 
of $z^n$ and $z^{*n}$ of the form 
$\sum_k\psi_{n,k}\dy$, $\psi_{n,k}\in\C$, and the functions 
$\sum_k\psi_{n,k}\delta_k(t)$ belong to $\FS$ since all sums 
are finite.                             
\phantom{   }                                      \hfill $\Box$    

A similar result can be obtained by considering the 
following set of (Borel measurable) functions
$$
 \Sy :=\{\,\psi : \sigma(\bar y)\rightarrow \C\,;\, 
 \sup_{t\in\sigma(\bar y)}|t^k\psi(t)|<\infty\ \,
  \mathrm{for\ all}\ k\in\Z\,\}. 
$$
\begin{thl}                         \label{SD}
The element 
$f=\sum_{n=0}^N z^n\psi_n(\bar y)+
     \sum_{n=1}^M\psi_{-n}(\bar y)z^{*n}$, 
$\psi_n \in \Sy$,  $N,M\in\N$,
belongs to $\BBA$. The operators $\psi(\bar y)$, 
$\psi\in\Sy$, satisfy on $D_\fA$ the commutation rules 
\begin{equation}                     \label{zpsi}
  z\psi(\bar y)=\psi(q^2\bar y)z,\quad 
  z^*\psi(\bar y)=\psi(q^{-2}\bar y)z^*.
\end{equation}
The linear space 
\begin{equation*}                     %\label{DefSD}
 \cS(D):=\left\{\,\sum_{n=0}^Nz^n\psi_n(\bar y)
  +\sum_{n=1}^M\psi_{-n}(\bar y)z^{*n}\,;\,\psi_k\in\Sy\ \,
  \mathrm{for\ all}  -M\leq k \leq N\,\right\}
\end{equation*}
forms a $\sut$-module *-subalgebra of $\ldop$.
\end{thl}
{\bf Proof.} 
By definition of $\BBA$, $a\psi b\in\BBA$ for all 
$a,b\in\fA$ whenever $\psi\in\BBA$. Fix $a\in\fA$. 
From the proof of  Lemma \ref{l-top}(i), we know that 
$\{||\cdot||_{\bar T^n}\}_{n\in\N}$, $T=1+y^2+y^{-2}$, 
generates the graph topology on $D_\fA$, so there exist 
$n_a\in \N$ such that $||a\varphi||\leq||T^{n_a}\varphi||$ 
for all $\varphi\in D_\fA$. Consequently, 
$||aT^{-n_a}\varphi||\leq||\varphi||$, 
%%%for all $\varphi\in D_\fA$, 
hence $\ov{aT^{-n_a}}$ and $\ov{T^{-n_a}a^*}$ are bounded. 
The operators $\ov{\py T^m}$, $\psi_n\in\Sy$, $m\in\N$, are bounded by 
the definition of $\Sy$, and $\bar T^{-1}$ is of trace class. 
From this facts, we conclude that 
$$
\ov{a\py b}=\ov{aT^{-n_a}}\, 
 \ov{\py T^{n_a+n_b+1}}  \bar T^{-1} \ov{T^{-n_b}b}
$$
is of trace class. This shows that the operator $f$ from
Lemma \ref{SD} belongs to $\BBA$. 

The commutation relations \rf[zpsi] are satisfied  
if we restrict the operators to $D\subset D_\fA$. 
Consider the O*-algebra generated by the elements 
$\psi(\bar y)\lceil D$, $\psi\in\Sy$, and 
$a\lceil D$, $a\in\fA$. Since the operators 
$\psi(\bar y)$ are bounded, the closure of this algebra 
is contained in $\ldop$. Taking the closure of an 
O*-algebra does not change the commutation relations, 
hence Equation \rf[zpsi] holds. 

Recall that $\fA$ is the linear span of operators
$z^np_n(y,y^{-1})$ and $p_{-n}(y,y^{-1})z^{*n}$, 
where $p_n(y,y^{-1})$ and $p_{-n}(y,y^{-1})$ are polynomials 
in $y$ and $y^{-1}$. 
Notice, furthermore, that $p(t,t^{-1})\psi(t)\in\Sy$ 
for all $\psi(t)\in\Sy$ and all polynomials $p(t,t^{-1})$. 
Now it follows from \rf[yz], \rf[zz*], \rf[zpsi], 
and the definition of $\cS(D)$ that $\cS(D)$ is stable under 
the $\sut$-action defined in Lemma \ref{discact}. 
Similarly, using \rf[zz*], \rf[zpsi], and 
the definition of $\cS(D)$, it is easy to check that $\cS(D)$ 
forms a *-algebra. Therefore, by Lemma \ref{discact}, 
$\cS(D)$ is a $\sut$-module *-algebra. 
                                          \hfill $\Box$
\mn

The description of $\FFD$ and $\cS(D)$ by functions 
$\psi:\sigma(\bar y)\rightarrow\C$ suggests that we can consider the 
elements of $\FFD$ and $\cS(D)$ as infinitely differentiable 
functions on the quantum disc with compact support and 
which are rapidly decreasing, respectively. 
Notice that $\FFD\neq\FFDA$ 
(e.g., $\eta\otimes\eta \notin \FFD$ for 
$\eta =\sum_{n=0}^\infty \alpha_n\eta_n\in\D_\fA$
if an infinite number of $\alpha_n$ are non-zero), 
and $\Sy\neq\BBA$ (e.g., $f=\sum_{k=0}^\infty 
\mathrm{exp}(-\bar y^{2^k})\dy z^{*k}\in \BBA$, 
$f\notin\cS(D)$). 

Clearly, $\FFD\subset\cS(D)$. On $\cS(D)$, the invariant integral 
can be expressed nicely in terms of the Jackson integral. 
The Jackson integral 
%on the interval $[0,1]$ and  $[0,\infty)$
is defined by 
$$
{\int^1_0}\varphi(t)d_qt=(1-q)
{\sum^{\infty}_{k=0}}\varphi(q^k)q^k\quad \mbox{{\rm and}}\quad
%$$
%and on $[0,\infty)$ by 
%$$
{\int^{\infty}_0}\varphi(t)d_qt=(1-q)
{\sum^{\infty}_{k=-\infty}}\varphi(q^k)q^k.
$$
\begin{thp}                                     \label{SDhaar}
Suppose that 
$\psi=\sum_{n=1}^Nz^n\psi_n(\bar y)+\psi_0(\bar y)
  +\sum_{n=1}^M\psi_{-n}(\bar y)z^{*n}\in\cS(D)$.
Let $h$ denote the invariant integral defined in 
Proposition \ref{haar}. For irreducible type $(I)$ 
rep\-re\-sen\-ta\-tions, we have 
$$
 h(\psi)= c (1-q^2)^{-1}{\int^1_0}\psi_0(t)t^{-2}d_{q^2}t,
$$
and, for irreducible type $(II)_\alpha$ rep\-re\-sen\-ta\-tions, we 
have
$$
h(\psi)= cq^{-2\alpha}(1-q^2)^{-1}
{\int^\infty_0}\psi_0(-q^{2\alpha}t)t^{-2}d_{q^2}t.
$$
\end{thp}
{\bf Proof.} Since 
$\ip{\eta_k}{z^n\psi_n(\bar y)y^{-1}\eta_k}=
\ip{\eta_k}{\psi_{-n}(\bar y)z^{*n}y^{-1}\eta_k}=0$ 
for all $n\neq 0$, we obtain 
\begin{align*}
h(\psi)
 &=c\, \tr \ov{\psi y^{-1}}=
c \sum_{k=0}^\infty
\ip{\eta_k}{\psi_0(\bar y)y^{-1}\eta_k}=
c{\sum^{\infty}_{k=0}}\psi_0(q^{2k})q^{-2k}\\
 &=c (1-q^2)^{-1}{\int^1_0}\psi_0(t)t^{-2}d_{q^2}t,
\end{align*}
for the type $(I)$ rep\-re\-sen\-ta\-tion and
\begin{align*}                  
h(\psi)&=c\, \tr \ov{\psi y^{-1}}=
c \sum_{k=-\infty}^\infty
\ip{\eta_k}{\psi_0(\bar y)y^{-1}\eta_k}=
c{\sum^{\infty}_{k=-\infty}}
\psi_0(-q^{2\alpha}q^{2k})q^{-2(\alpha+k)} \\
 &=cq^{-2\alpha}(1-q^2)^{-1}
{\int^\infty_0}\psi_0(-q^{2\alpha}t)t^{-2}d_{q^2}t
\end{align*}
for type $(II)_\alpha$ rep\-re\-sen\-ta\-tions.      \hfill $\Box$
%
%
%---------------------------------------------------------------------
\subsection{Application: differential calculus}
%---------------------------------------------------------------------  
%
                                                  \label{calc}
The bimodule structure of a first order differential *-calculus 
$(\Gamma, \dd)$ over $\qd$ has been described 
in \cite{VK1} and \cite{S-com}.  
The commutation relations are given by 
\begin{equation*}                        
 \dd z\,z=q^2z\,\dd z,\quad \dd z\,z^*=q^{-2}z^*\,\dd z,\quad
 \dd z^*\,z=q^2z\,\dd z^*,\quad \dd z^*\,z^*=q^{-2}z^*\,\dd z^*.
\end{equation*}
Our aim is to extend this FODC to the classes of integrable functions 
on the quantum disc defined in Subsection \ref{topasp}. 
To this end, we use a commutator rep\-re\-sen\-ta\-tion of the FODC. 
A faithful commutator rep\-re\-sen\-ta\-tion of the above 
differential calculus can be found in \cite{S-com} and is 
obtained as follows. Given a *-rep\-re\-sen\-ta\-tion $\pi$ of $\qd$ 
from Subsection \ref{topasp}, consider the direct sum 
$\rho:=\pi\oplus\pi$ on $D\oplus D\subset\Hh\oplus\Hh$ and 
set 
\[ C:=(1-q^2)^{-1} \left(
      \begin{array}{cc}    0     &  \pi(z) \\
                      \pi( z^*)  &    0    \end{array} \right).
\] 
Then the 
differential mapping $\dd_{\rho,C}$ defined in \rf[dop]
is given by 
\[ 
 \dd_{\rho,C}(f)=\im [C,\rho(f)]=(1-q^2)^{-1}\im \left(
 \begin{array}{cc}    0     &  \pi(zf-fz) \\
          \pi( z^*f-fz^*)   &    0    \end{array} \right),\ f\in\qd. 
\]
Clearly, $C\in\ldop[D\oplus D]$, so we can extend 
$\dd_{\rho,C}$ to $\ldop[D\oplus D]$, that is, 
$$
 \dd_{\rho,C}(x):=\im [C,x],\quad x\in \ldop[D\oplus D].
$$
The same formula applies to any *-subalgebra of $\ldop[D\oplus D]$. 
Notice that we can consider $\ld$ as a *-subalgebra of
$\ldop[D\oplus D]$ by identifying $A\in\ld$ with the 
operator $A\oplus A$ acting on $D\oplus D$. In particular, 
the algebras $\FFD$ and $\BBA$ from Proposition \ref{haar} 
become *-subalgebras of $\ldop[D\oplus D]$. 
In this way, we obtain a FODC over these algebras.  

For $z$ and $z^*$, we have 
\[
\dd_{\rho,C}(z)=\im \left(
 \begin{array}{cc}    0     &    0 \\
          \pi( y)   &    0    \end{array} \right),\quad
\dd_{\rho,C}(z^*)=\im \left(
 \begin{array}{cc}  0   &    -\pi(y)    \\
                    0   &    0    \end{array} \right).
\]
For functions $\psi(\bar y)$, the differential mapping 
$\dd_{\rho,C}$ can be expressed in terms of the 
$q$-differential operator $\Dq$ defined by 
$\Dq f(x)=(x-qx)^{-1}(f(x)-f(qx))$. It follows from 
$$
(1-q^2)^{-1}(z\psi(\bar y)-\psi(\bar y)z)=
zy(y-q^2y)^{-1}(\psi(\bar y)-\psi(q^2\bar y))
=z\Dqq\psi(\bar y)y,
$$
$$
(1-q^2)^{-1}(z^*\psi(\bar y)-\psi(\bar y)z^*)=
y(y-q^2y)^{-1}(\psi(q^2\bar y)-\psi(\bar y))z^*
=-q^{-2}\Dqq\psi(\bar y)z^*y
$$
that
$$
\dd_{\rho,C}(\psi(\bar y))=
-\im \rho(z)\Dqq\psi(\bar y)\dd_{\rho,C}(z)
-\im q^{-2}\Dqq\psi(\bar y)\rho(z^*)\dd_{\rho,C}(z^*).
$$
In particular, the ``$\delta$-distributions'' $\dy$ 
are differentiable.

%
%++++++++++++++++++++++++++++++++++++++++++++++++++++++++++++++++++++++++
\section{Quantum $(n,1)$-matrix ball}
%++++++++++++++++++++++++++++++++++++++++++++++++++++++++++++++++++++++++
                                                     \label{balls}
%
%------------------------------------------------------------------------
\subsection{Algebraic relations}
%------------------------------------------------------------------------
                                                         \label{sec-ball}
Let $n\in\N$ and $q\in(0,1)$. 
We denote by $\qm$ the *-algebra generated by 
$z_1,\ldots,z_n,z^*_1,\ldots,z^*_n$  obeying the relations 
\begin{eqnarray}                       
  z_k z_l &=& q z_l z_k,\quad k<l,       \label{ball1} \\
  z^*_l z_k &=& q z_k z^*_l,\quad k\neq l, \label{ball2}\\
  z^*_k z_k &=& q^2 z_k z^*_k                 \label{ball3} 
  -(1-q^2)\sum^n_{j=k+1}z_j z^*_j+(1-q^2),\quad k<n,  \\
 z^*_n z_n &=& q^2 z_n z^*_n + (1-q^2).     \label{ball4} 
\end{eqnarray}  

Equations \rf[ball1]--\rf[ball4] are called 
{\it twisted canonical commutation relations} \cite{PW} 
and $\qm$ is also known as {\it $q$-Weyl algebra} \cite{KS}. 
Here we consider it as a special case of the 
quantum matrix balls introduced in \cite{VK2} because the 
$\sutn$-action on $\qm$ defined below is taken from the latter.  

The following hermitian elements $Q_k$ will play a
crucial role throughout this section. Set 
\begin{equation}                              \label{Q}
   Q_k:= 1-\sum_{j=k}^nz_jz_j^*,\,\ k\leq n,\quad \, 
   Q_{n+1}:=1. 
\end{equation}
Equations \rf[ball3], \rf[ball4], and \rf[Q] imply immediately 
\begin{equation}                              \label{zzQ}
 z_k^*z_k-q^2z_kz_k^* =(1-q^2)Q_{k+1},\quad \,
     z_k^*z_k-z_kz_k^* =(1-q^2)Q_{k}. 
\end{equation}
Taking the difference of the first  with the second and of the first 
with $q^2$ times the second equation gives
\begin{equation}                              \label{zQQ} 
     z_kz_k^*=Q_{k+1}-Q_k,\quad\, z_k^*z_k=Q_{k+1}-q^2Q_k.
\end{equation}
Furthermore, one easily shows by using 
Equations \rf[ball1]--\rf[Q] that
\begin{equation}                              \label{Qz}
 Q_kz_j=z_jQ_k,\,\ j<k,\quad \, Q_kz_j=q^2z_jQ_k,\,\ j\geq k,
\end{equation}
\begin{equation}                              \label{Qz*}
Q_kz_j^*=z_j^*Q_k,\,\ j<k,\quad\,Q_kz_j^*=q^{-2}z_j^*Q_k,\,\ j\geq k.
\end{equation}
As a consequence, 
\begin{equation}                              \label{QQ}
              Q_kQ_l=Q_lQ_k ,\,\quad \,
\mathrm{for\ all}\ \, k,l\leq n+1.
\end{equation}

For $I=(i_1,\ldots,i_n)\in\N_0^n,\ J=(j_1,\ldots,j_n)\in\N_0^n$, 
set                                                    \label{index}
$z^I:=z_1^{i_1}\cdots z_n^{i_n}$,\ \,$z^{*J}:=z_1^{*j_1}\cdots z_n^{*j_n}$ 
and define                                  
$I\cdot J=(i_1j_1,\ldots,i_nj_n)\in\N_0^n$. We write $0$ instead of 
$(0,\ldots,0)$.
It follows from \rf[Qz]--\rf[QQ] together with the defining 
relations \rf[ball1]--\rf[ball4] that each $f\in\qm$ can be expressed  
as a finite sum 
\begin{equation}                                \label{f}
f=\sum_{I\cdot J=0}z^Ip_{IJ}(Q_1,\ldots,Q_n)z^{*J}
\end{equation}
with polynomials $p_{IJ}(Q_1,\ldots,Q_n)$ 
in $Q_1,\ldots,Q_n$.

The $\sutn$-action $\ang$ on $\qm$ which turns $\qm$ into a 
$\sutn$-module *-algebra is given by the 
following formulas \cite{VK2}. 
\begin{align*}
j&\neq n: & E_j\ang z_{j+1} &\ \,=\ \,q^{-1/2}z_j,                  & 
                          E_j\ang z_k &\ \,=\ \,0,\quad k\neq j+1,  \\
& & E_j\ang z_{j}^* &\ \,=\ \,-q^{-3/2}z_{j+1}^*,                   & 
                          E_j\ang z_k^* &\ \,=\ \,0,\quad k\neq j,  \\
& & F_j\ang z_{j} &\ \,=\ \,q^{1/2}z_{j+1},                         &
                          F_j\ang z_k &\ \,=\ \,0,\quad k\neq j,    \\ 
& & F_j\ang z^*_{j+1} &\ \,=\ \,-q^{3/2}z^*_{j},                    &
                          F_j\ang z^*_k &\ \,=\ \,0,\quad k\neq j+1, \\
\lefteqn{
K_j\ang z_j\ \,=\ \,qz_j,\quad 
K_j\ang z_{j+1} \ \,=\ \,q^{-1}z_{j+1}, \,\quad
K_j\ang z_k\ \,=\ \,z_k,\,\ k\neq j,j+1,  }\\
\lefteqn{ 
 K_j\ang z_j^*\ \,=\ \,q^{-1}z_j^*,\quad 
K_j\ang z_{j+1}^* \ \,=\ \,qz_{j+1}^*, \,\quad 
K_j\ang z_k^*\ \,=\ \,z_k^*,\,\ k\neq j,j+1, } 
\end{align*}
\begin{align*}
j&=n: & E_n\ang z_n &\ \,=\ \,-q^{1/2}z_n^2,                           &
     k&<n:&             E_n\ang z_k &\ \,=\ \,-q^{1/2}z_nz_k,       \\
& & E_n\ang z^*_n &\ \,=\ \,q^{-3/2},                                  & 
  &&               E_n\ang z^*_k &\ \,=\ \,0,                 \\
& & F_n\ang z_n &\ \,=\ \,q^{1/2},                                     &
   &&              F_n\ang z_k &\ \,=\ \,0,                    \\
& & F_n\ang z^*_n &\ \,=\ \,-q^{5/2}z^{*2}_n                           &
   &&              F_n\ang z^*_k &\ \,=\ \,-q^{5/2}z^*_kz^*_n,\qquad  \\ 
& & K_n\ang z_n &\ \,=\ \,q^2z_n,                                      &
    &&             K_n\ang z_k &\ \,=\ \,qz_k,               \\
& & K_n\ang z^*_n &\ \,=\ \,q^{-2}z^*_n,                               &
      &&           K_n\ang z^*_k &\ \,=\ \,q^{-1}z^*_k.
\end{align*}
If $n=1$, we recover the relations of the quantum disc. 
For $n>1$, we obtain by omitting  
the elements $K_n$, $K_n^{-1}$, $E_n$, and $F_n$ 
a $\sutm$-action on $\qm$  such that 
$\qm$ becomes a $\sutm$-module *-algebra.
%In the following, we assume that $n>1$.  
Notice that, by Equation \rf[modalg] and \rf[delta], it is 
sufficient to describe the action on generators. 
%
%
%
%++++++++++++++++++++++++++++++++++++++++++++++++++++++++++++++++++++++++
\subsection{Representations of the *-algebra $\qm$}
%++++++++++++++++++++++++++++++++++++++++++++++++++++++++++++++++++++++++
%
Irreducible *-rep\-re\-sen\-ta\-tions of the twisted canonical 
commutation relations have been  
classified in 
%%%%%%%%%%\cite{OS} and 
\cite{PW} %%%%%%%. The latter includes the 
%%%%%%%%classification of the irreducible 
%%%%%%%%unbounded *-rep\-re\-sen\-ta\-tions 
under the 
condition that $1-Q_1$ is essentially self-adjoint. 
In this subsection, we study admissible 
*-rep\-re\-sen\-ta\-tions of  the twisted canonical commutation relations
without requiring the rep\-re\-sen\-ta\-tion to be irreducible. 

Remind our notational conventions from Subsection \ref{oppre} regarding 
direct sums of a Hilbert space $\K$. 
Let $A$ be a self-adjoint operator on $\K$ such that 
$\sigma(A)\sqsubseteq  (q^2,1]$. 
Then the expression $\mu_{j}(A)$, $j\in\Z$, stands for the operator 
$\mu_{j}(A)=(1+q^{-2j}A)^{1/2}$. 
We shall also abbreviate $\lambda_j=(1-q^{2j})^{1/2}$ and 
$\beta_j=(q^{-2j}-1)^{1/2}$ for $j\in\N_0$. 
\begin{thp}                                       \label{DarMat}
Assume that $m,k,l\in\N_0 $ such that $m+l+k=n$. 
Let $\K$ denote a Hilbert space.
Set
$$
\Hh:=\oplus^\infty_{i_n,\ldots, i_{n-m+1}=0}\oplus^\infty_{i_{k}=-\infty}
     \oplus^\infty_{i_{k-1},\ldots,i_1=1}\Hh_{i_n\ldots i_1}, 
$$
where $\Hh_{i_n\ldots\ i_1}=\K$, and  
$$
D:=\Lin \{\eta_{i_n\ldots i_1}\,;\,\eta\in\K,\ 
i_n,\ldots,i_{n-m+1}\in\N_0,\ i_{k}\in\Z,\ i_{k-1},\ldots, i_1\in\N \}.
$$ 
(For $l>0$, we retain the notation $\eta_{i_n\ldots i_1}$ and do not write 
$\eta_{i_n\ldots i_{n-m+1},i_{k}\ldots i_1}$.)
Consider the operators $z_1,\ldots,z_n$ acting 
on $D$ by 
\begin{align*}
\lefteqn{(m,0,k):}\\ 
&& z_j\eta_{i_n\ldots i_1} &= q^{i_{j+1}+\ldots+i_1}
\lambda_{i_j+1}\eta_{i_n\ldots i_j+1\ldots i_1},\,\quad  
\mbox{if}\ \,k<j\le n,\\
&& z_{k}\eta_{i_n\ldots i_1} &= q^{i_{k+1}+\ldots+i_n}
\mu_{i_{k}-1}(A^2)\eta_{i_n\ldots i_{k}-1\ldots i_1},\\
&& z_j\eta_{i_n\ldots i_1} &= q^{-(i_{j+1}+\ldots+i_{k})+
(i_{k+1}+\ldots+i_n)}\beta_{i_j-1} A\eta_{i_n\ldots i_j-1\ldots i_1},
\,\quad  \mbox{if}\  \, 1\le j<k, 
\hspace{18.3pt}
\end{align*}
and,\, for\, $l>0$, 
\begin{align*}
\lefteqn{(m,l,k):}\\ 
&& z_j\eta_{i_n\ldots i_1} &= q^{i_{j+1}+\ldots+i_1}
\lambda_{i_j+1}\eta_{i_n\ldots i_j+1\ldots i_1},\,\quad 
\mbox{if}\ \,n-m<j\le n,\\
&& z_{n-m}\eta_{i_n\ldots i_1} &= 
q^{i_{n-m+1}+\ldots+i_n}v\eta_{ i_n\ldots i_{k}-1\ldots i_1},\\       
&& z_j &\equiv 0,\,\quad\mbox{if}\ \,k<j<n-m,\\
&& z_{k}\eta_{i_n\ldots i_1} &= q^{-i_{k}+i_{n-m}+\ldots+i_n}
 A\eta_{i_n\ldots i_{k}-1\ldots i_1},\\
&& z_j\eta_{i_n\ldots i_1} &= q^{-(i_{j+1}+\ldots+i_{k})+
(i_{n-m+1}+\ldots i_n)}\beta_{i_j-1} A\eta_{i_n\ldots i_j-1\ldots i_1},
\quad  \mbox{if}\  \, 1\le j<k.
\end{align*}
(If $k=0$, then the indices $i_1,\ldots, i_{k}$ are omitted;
similarly, if $m=0$, then the indices $ i_{n-m+1},\ldots,i_n$ 
are omitted.)
In both series, $A$ denotes a self-adjoint operator acting on the Hilbert space
$\K$ such that $\sigma(A)\sqsubseteq  (q,1]$.
In the series $(m,l,k),\ l>0$, $v$ is a unitary operator on $\K$
such that $Av=vA$. 

Then the operators $z_1,\ldots,z_n$ define 
a *-rep\-re\-sen\-ta\-tion
of $\qm$, where the action of $z_j^*$, $j=1,\ldots,n$ is obtained
by restricting the adjoint of $z_j$ to $D$. 
Representations belonging to different series $(m,k,l)$ or to different
operators $A$ and $v$ are not unitarily equivalent.
A rep\-re\-sen\-ta\-tion of this series is irreducible 
if and only if ${\K=\C}$.
In this case, $v$ is a complex number of modulus one and
$A\in (q,1]$. Only the rep\-re\-sen\-ta\-tions $(m,0,k)$ are faithful.
\end{thp}

\noindent
{\bf Proof.} Direct calculations show that  the formulas given in 
Proposition \ref{DarMat} define a 
*-rep\-re\-sen\-ta\-tion of $\qm$. Clearly, if a *-representation of these 
series is irreducible, then $A$ and $v$ must be complex numbers 
and $\K=\C$. The converse statement was shown in \cite{PW}. 
That the rep\-re\-sen\-ta\-tions $(m,0,k)$ are 
faithful is proved by showing that
for each ${x\in\qm}$, $x\neq0$, there exist $\eta_{i_n\ldots i_1},
\eta_{j_n\ldots j_1}\in \Hh$ such that the matrix element
$\langle \eta_{i_n\ldots i_1},x\eta_{j_n\ldots j_1}\rangle$ is non-zero. 
The vectors can easily be found by writing $x$ in the standard 
form (\ref{f}) and observing that $z_j$, $z_j^*$ act as shift operators. 
We omit the details. The other assertions of the proposition are obvious.
    \mbox{ } \mbox{ }  \mbox{ } \mbox{ }    \mbox{ }     
                                    \hfill $\Box$
\mn\\
{\bf Remarks.}
The operators $Q_j$ are given by
\begin{align}                                               
& Q_j\eta_{i_n\ldots i_1}=q^{2(i_j+\ldots+i_n)}\eta_{i_n\ldots i_1},\,
\quad \mbox{if}\  \, n-m<j\le n,                         \label{y1}\\
& Q_j\equiv 0,\,\quad \mbox{if}\  \,  k<j\le n-m,        \label{y2}\\
& Q_j\eta_{i_n\ldots i_1}=-q^{-2(i_j+\ldots+i_{k})+2(i_{n-m+1}+\ldots+i_n)}
A^2\eta_{i_n\ldots i_1},\,\quad \mbox{if}\  \,  1\le j\le k.  \label{y3}
\end{align}
The numbers $m,l,k\in\N_0$ correspond to the signs of the operators $Q_j$,
that is, we have $Q_n\ge \ldots\ge Q_{n-m+1}>0$ if $m>0$,
$ Q_{n-m}=\ldots=Q_{k+1}=0$ if $l>0$, and
$0>Q_{k}\ge\ldots \ge Q_{1}$ if $k>0$.
The only bounded rep\-re\-sen\-ta\-tions are the series 
$(m,l,0)$.
\mn

We now give a constructive method for finding ``admissible'' 
*-representations of $\qm$. 
In view of \rf[Qz]--\rf[QQ], 
the assumptions on  {\it admissible} *-rep\-re\-sen\-ta\-tions of 
the *-al\-ge\-bra $\qm$ 
will include the following two conditions: 
First, the closures of the operators
$Q_k$, $k=1,\ldots,n$, are  
self-adjoint and strongly commute. Second, 
$\varphi(\bar Q_k)z_j\subset z_j\varphi(\bar Q_k)$, $j<k$,
and $\varphi(\bar Q_k)z_j\subset z_j\varphi(q^2\bar Q_k)$, $j\geq k$,
for all complex functions $\varphi$ which are measurable with respect to 
the spectral measure of $\bar Q_k$ and which have  at most 
polynomial growth. 
In the course of the argumentation, 
we shall impose further regularity conditions on the operators. 
The outcome will precisely be the series 
of Proposition \ref{DarMat}. So, if one takes as admissible 
*-rep\-re\-sen\-ta\-tions those which satisfy all regularity conditions, 
then Proposition \ref{DarMat} states that any admissible 
*-rep\-re\-sen\-ta\-tion of $\qm$ is a direct sum of 
*-rep\-re\-sen\-ta\-tion which are determined by the formulas 
of the series $(m,l,k)$, $m+l+k=n$. 
The argumentation is based on a reduction procedure.  

Observe that $z_n$ satisfies the relation of the quantum
disc $\cO_q(\bar{U})$. 
The ``admissible'' rep\-re\-sen\-ta\-tions of this relation 
are given by Lemma \ref{L3} 
and correspond to the cases
$(1,0,0)$, $(0,1,0)$, and $(0,0,1)$. 

Now let $0<m<n$.  
Suppose that we are given a *-rep\-re\-sen\-ta\-tion of $\qm$
such that the operators
$z_n,\ldots, z_{n-m+1}$
act on 
%%%%the Hilbert space 
$\Hh=\oplus^\infty_{i_n,\ldots,i_{n-m+1}=0}\Hh_{i_n\ldots i_{n-m+1}}$
%%%%%%%%%%%%%, $\Hh_{i_n\ldots i_{n-m+1}}=\Hh_{0\ldots0}$, 
by the formulas of the series $(m,0,0)$,
where all $\Hh_{i_n\ldots i_{n-m+1}}$ are equal to a given Hilbert space, 
say $\Hh_{0\ldots0}$. 
Fix $f_n,\ldots,f_{n-m+1}\in \N_0$.
Since the rep\-re\-sen\-ta\-tion is assumed to be admissible, it follows from
$z_{n-m}Q_j=Q_j z_{n-m}$, $n-m<j\le n$, 
and \rf[y1]
that the operator 
$ z_{n-m}$ maps the Hilbert spaces
\begin{align*}
&\Hh(f_n):=\Lin\{ \eta_{i_n\ldots i_1}\in\Hh\,;\, i_n=f_n\}, \\
&\Hh(f_n,f_{n-1}):=\Lin\{ \eta_{i_n\ldots i_1}\in\Hh\,;\, 
i_n=f_n, \ i_n+i_{n-1}=f_{n-1}\},\ldots,     \\                  
&\Hh(f_n,\!...,f_{n-m+1})
:=\Lin\{ \eta_{i_n\ldots i_1}\in\Hh\,;\, 
i_n\!=\!f_n,\ldots,  i_n+...+i_{n-m+1}\!=\!f_{n-m+1}\}
\end{align*}
%
%$\oplus_{i_n=f_n}\oplus^\infty_{i_{n-1},\ldots,i_{n-m+1}=0}
%\Hh_{i_n\ldots i_{n-m+1}}$, $\ldots\ $, $\oplus_{i_n=f_n}{\ldots}
%\oplus_{i_n+\ldots+i_{n-m+1}=f_{n-m+1}} \Hh_{i_n\ldots i_{n-m+1}}$
into itself.
But the $m$ equations $i_n=f_n$, $\ldots\ $, $i_n+\ldots+i_{n-m+1}=f_{n-m+1}$
determine uniquely the numbers $i_n,\ldots,i_{n-m+1}$, therefore $z_{n-m}$
maps each $\Hh_{i_n\ldots i_{n-m+1}}$ into itself.
Write 
$$
z_{n-m}\eta_{i_n\ldots i_{n-m+1}}=
Z_{i_n\ldots i_{n-m+1}}\eta_{i_n\ldots i_{n-m+1}}
$$
with operators $Z_{i_n\ldots i_{n-m+1}}$ 
acting on $\Hh_{0\ldots0}$. 
Applying $ z_{n-m} z_j=q z_j z_{n-m}$, $n-m<j\le n$, 
to vectors $\eta_{i_n\ldots
i_{n-m+1}}$ gives
\begin{eqnarray*}
\lefteqn{
q^{i_{j+1}+\ldots+i_n}
\lambda_{i_j+1}Z_{i_n\ldots i_j+1\ldots i_{n-m+1}}
\eta_{i_n\ldots i_j+1\ldots i_{n-m+1}}= }\\
&& \qquad\qquad\qquad\qquad\qquad
q q^{i_{j+1}+\ldots+i_n}
\lambda_{i_j+1}Z_{i_n\ldots i_j\ldots i_{n-m+1}}
\eta_{i_n\ldots i_j+1\ldots i_{n-m+1}},
\end{eqnarray*}
hence $ Z_{i_n\ldots i_j+1 \ldots i_{n-m+1}}=
q Z_{i_n\ldots i_j \ldots i_{n-m+1}}$.
From this, we conclude
$$ 
Z_{i_n\ldots i_{n-m+1}}= q^{i_{n-m+1}+\ldots+i_n}Z_{0\ldots0}.
$$
On $\Hh_{0\ldots0}$, the relation
$ z_{n-m}^*z_{n-m}-q^2 z_{n-m} z_{n-m}^* =(1-q^2)Q_{n-m+1}$ gives
$$ 
Z_{0\ldots0}^* Z_{0\ldots0}-q^2 Z_{0\ldots0}Z_{0\ldots0}^*=1-q^2.
$$
Here and subsequently, we suppose that operators satisfying 
this relation also satisfy the assumptions of Lemma \ref{L3}. 
By a slight reformulation of Lemma \ref{L3}, 
we get the following  three series of *-rep\-re\-sen\-ta\-tions of the
last equation:
\begin{enumerate}
\item
$\Hh_{0\ldots 0}= \oplus^\infty_{i_{n-m}=0} \Hh_{0\ldots 0,i_{n-m}},\
\Hh_{0\ldots 0,i_{n-m}}=\Hh_{0\ldots0,0},$
$$
Z_{0\ldots0}\zeta_{i_{n-m}}=\lambda_{i_{n-m}+1}\zeta_{i_{n-m}+1};
$$
\item
$ Z_{0\ldots 0}= v$, where $v$ is a unitary operator on
$\Hh_{0\ldots0}$;
\item
$\Hh_{0\ldots 0}=\oplus^\infty_{i_{n-m}=-\infty} \Hh_{0\ldots 0,i_{n-m}},\
\Hh_{0\ldots 0,i_{n-m}}= \Hh_{0\ldots 0,0},$
$$
Z_{0\ldots 0}\zeta_{i_{n-m}}=
\mu_{i_{n-m}-1}(A^2)\zeta_{i_{n-m}-1},
$$
where $A$ is a self-adjoint operator on $\Hh_{0\ldots 0,0}$ such
that $\sigma(A)\sqsubseteq (q^2,1]$.
\end{enumerate}
Inserting these formulas into the rep\-re\-sen\-ta\-tion $(m,0,0)$ shows that
the cases (i), (ii), and (iii) correspond to 
a rep\-re\-sen\-ta\-tion of the operators 
$ z_n,\ldots, z_{n-m}$ of the series
$(m+1,0,0)$, $(m,1,0)$, and $(m,0,1)$, respectively.

Next, let $ m,l\in\N_0$ such that $ m+l<n$ and $l>0$. Set $k=n-m-l$. 
Suppose that the operators $ z_n,\ldots, z_{k+1}$ act on the Hilbert space
$\Hh=\oplus^\infty_{i_n,\ldots,i_{n-m+1}=0}\Hh_{i_n\ldots i_{n-m+1}}$
by the formulas of the series $(m,l,0)$.
As in the case $(m,0,0)$, we conclude from
$ z_{k} Q_j = Q_j z_{k}$, $k<j\le n$, 
that $z_{k}$ acts on $\Hh_{i_n\ldots i_{n-m+1}}$ by
$$ 
z_{k}\eta_{i_n\ldots i_{n-m+1}}= q^{i_{n-m+1}+\ldots+i_n}Z_{0\ldots 0}
\eta_{i_n\ldots i_{n-m+1}},
$$
where $Z_{0\ldots0}$ is an operator acting on $\Hh_{0\ldots0}$. 
As $ Q_{k}=Q_{k+1}-z_{k}z_{k}^*$ and $Q_{k+1}=0$ by \rf[zQQ] 
and \rf[y2], we have
$ Q_{k}\le 0 $. 
The assumptions on admissible *-rep\-re\-sen\-ta\-tions imply that 
${\rm  ker}~Q_{k}$ is reducing,
Thus we can consider the cases $Q_{k}=0 $ and $Q_{k}<0$ separately.  

First let $Q_{k}=0$.
Then $ z_{k}z_{k}^*=Q_{k+1}-Q_{k}=0$,
hence $z_{k}= z_{k}^*=0$. Consequently, 
we have obtained a rep\-re\-sen\-ta\-tion
of the type $(m,l+1,0)$. 

Now  assume that $Q_{k}<0$. Then, by \rf[zzQ], 
$z_{k}^* z_{k}-q^2 z_{k} z_{k}^*=0$ and
$z_{k}^* z_{k}- z_{k} z_{k}^*= (1-q^2)Q_{k}$. 
Inserting the first equation into the second one gives 
$z_{k} z_{k}^*=-Q_{k}$, hence
$|z_{k}^*|=|Q_{k}|^{1/2}$.
Evaluating $z_{k}^* z_{k}-q^2 z_{k}z_{k}^*$ on
$\Hh_{0\ldots 0} $ yields $Z_{0\ldots 0}^* Z_{0\ldots 0}-q^2
Z_{0\ldots 0} Z_{0\ldots 0}^*=0$.
Let $Z_{0\ldots 0}= U|Z_{0\ldots 0}|$ be the polar decomposition of the
closed operator $Z_{0\ldots 0}$. Since 
${\rm ker}~z_{k}^*={\rm ker}~z_{k} z_{k}^*
={\rm ker}~Q_{k}=\{0\}$ and
${\rm ker}~z_{k}={\rm ker}~z_{k}^*z_{k}
={\rm ker}~z_{k}z_{k}^*=\{0\}$, $U$ is unitary.
From $Z_{0\ldots 0}^*Z_ {0\ldots 0}-q^2Z_{0\ldots 0}Z_{0\ldots 0}^*=0$, 
it follows that $|Z_{0\ldots 0}|^2 =q^2U|Z_{0\ldots 0}|^2U^*$
and therefore $|Z_{0\ldots 0}|U =qU|Z_{0\ldots 0}|$.
The rep\-re\-sen\-ta\-tions of this relation are given by Lemma \ref{L2}.
It states that
$\Hh_{0\ldots 0}=\oplus^\infty_{i_{k}=-\infty}\Hh_{0\ldots 0,i_{k}}$,
$\Hh_{0\ldots 0,i_{k}}=\Hh_{0\ldots0,0}$, and the operators
$U$ and $|Z_{0\ldots 0}|$ act as
$$
U\zeta_{i_{k}}=\zeta_{i_{k}-1},\quad \,
|Z_{0\ldots 0}|\zeta_{i_{k}}=q^{-i_{k}}A\zeta_{i_{k}},
$$
where $A$ is a self-adjoint operator on
$\Hh_{0\ldots 0,0}$ such that $\sigma (A)\sqsubseteq(q,1]$.

We have not yet considered the relations
$z_{k}z_{n-m}=q z_{n-m} z_{k}$ and
$z_{k}^*z_{n-m}=q z_{n-m}z_{k}^*$.
On $\Hh_{0\ldots 0}$, this leads to
$Z_{0\ldots 0}v=qvZ_{0\ldots 0}$ and
$Z_{0\ldots 0}^*v=qvZ_{0\ldots 0}^*$.
Thus $Z_{0\ldots 0}^*Z_{0\ldots0}v=q^2vZ_{0\ldots0}^*Z_{0\ldots0}$,
or, since $v$ is unitary,
$|Z_{0\ldots0}|v=qv|Z_{0\ldots0}|$. This implies
$U^*v|Z_{0\ldots0}|=|Z_{0\ldots0}|U^*v$, hence the unitary operator
$U^*v$ leaves each space $\Hh_{0\ldots0,i_{k}}$ invariant.
Therefore there exist unitary operators $v_{i_{k}}$ on
$\Hh_{0\ldots0,i_{k}}$ such that
$U^*v\zeta_{i_{k}}=v_{i_{k}}\zeta_{i_{k}}$, hence
$v\zeta_{i_{k}}=v_{i_{k}}\zeta_{i_{k}-1}$.
From 
$$vU|Z_{0\ldots0}|=vZ_{0\ldots0}=q^{-1}Z_{0\ldots0}v
=q^{-1}U|Z_{0\ldots0}|v=Uv|Z_{0\ldots0}|,
$$ 
we conclude
$vU=Uv$ since ${\rm ker}~|Z_{0\ldots0}|=\{0\}$.
This implies 
$$
v_{i_{k}}\zeta_{i_{k}}=Uv\zeta_{i_{k}}
=vU\zeta_{i_{k}}=v_{i_{k}+1}\zeta_{i_{k}},
$$
thus $v_{i_{k}+1}= v_{i_{k}}$, consequently $v_{i_{k}}=v_0$
for all $i_{k}\in \Z$. Evaluating $|Z_{0\ldots0}|v=qv|Z_{0\ldots0}|$
on vectors $\zeta_{i_{k}}$ shows $v_0A=Av_0$. This determines the actions
of $z_n,\ldots, z_{k}$ completely.
Comparing the result with the action of the operators $ z_n,\ldots, z_{k}$
from the proposition shows that we have obtained 
a rep\-re\-sen\-ta\-tion of the type 
$(m,l,1)$.

We finally turn to a rep\-re\-sen\-ta\-tion of the type $(m,l,k)$, where
$m+l+k<n$ and $ k>0$. Set $s:=n-(m+l+k)$.
Suppose that the operators $ z_n,\ldots,z_{s+1}$ act on a Hilbert space
$\Hh=\oplus^\infty_{i_n,\ldots,i_{n-m+1}=0}\oplus^\infty_{i_{k}=-\infty}
     \oplus^\infty_{i_{k-1},\ldots, i_{s+1}=1}\Hh_{i_n\ldots i_{s+1}}$
by the formulas given in the proposition. Similarly to the case $(m,0,0)$, 
we conclude from
$z_{s}Q_j=Q_j z_{s}$, $s<j\le n$, that
$z_{s}$ maps 
each $\Hh_{i_n\ldots i_{s+1}}$ into itself. Write
$$
z_{s}\eta_{i_n\ldots i_{s+1}}=Z_{i_n\ldots i_{s+1}}\eta_{i_n\ldots
i_{s+1}}.
$$
On applying $ z_{s}Q_{s+1}=Q_{s+1}z_{s}$ to a vector
$\eta_{i_n\ldots i_{s+1}}$, it follows from Equation \rf[y3] that
$Z_{i_n\ldots i_{s+1}} A^2=A^2Z_{i_n\ldots i_{s+1}}$. 
Thus 
$Z_{i_n\ldots i_{s+1}} \mu_{i_{n-m}}(A^2)
=\mu_{i_{n-m}}(A^2)Z_{i_n\ldots i_{s+1}}$
because the
rep\-re\-sen\-ta\-tion is assumed to be admissible. 
By using this relations and 
evaluating
$z_{s}z_j=q z_j z_{s}$, $j>s$, on vectors
$\eta_{i_n\ldots i_{s+1}}$, we see that
$Z_{i_n\ldots i_j-1\ldots i_{s+1}}=q Z_{i_n\ldots i_j\ldots i_{s+1}}$
%$ s<j\le n-m$, 
and 
$Z_{i_n\ldots i_j+1\ldots i_{s+1}}=q Z_{i_n\ldots i_j\ldots i_{s+1}}$
if $ s<j\le k$ and  $ n-m<j\le n$, respectively.
Thus we can write
$$
Z_{i_n\ldots i_{s+1}}=q^{-(i_{s+1}+\ldots+i_{k})+(i_{n-m+1}+\ldots+i_n)}
Z_{0\ldots0}.
$$ 
When $l\neq0$, this formula together with
$z_{s} z_{n-m}=q z_{n-m}z_{s}$ implies
$Z_{0\ldots0}v=vZ_{0\ldots0}$.

The relation $z_{s}^* z_{s}-q^2z_{s}z_{s}^*=(1-q^2)Q_{s+1}$
reads on $\Hh_{0\ldots0}$ as
$$
Z_{0\ldots0}^*Z_{0\ldots0}-q^2Z_{0\ldots0}Z_{0\ldots0}^*=-(1-q^2)A^2.
$$
Setting $X_{0}:=Z_{0\ldots0}A^{-1}= A^{-1}Z_{0\ldots0}$ and replacing
$Z_{0\ldots0}$ by $X_{0}A$, we obtain
$$
X_{0}^*X_{0}-q^2X_{0}X_{0}^*=-(1-q^2).
$$
By Lemma \ref{L3}, 
the admissible *-representations of 
this relation can be described in the following way: 
$\Hh_{0\ldots0}=\oplus^\infty_{i_{s}=1}\Hh_{0\ldots0,i_{s}},\
\Hh_{0\ldots0,i_{s}}=\Hh_{0\ldots0,1}$, and $X_{0}$ acts by
$$
X_{0}\zeta_{i_{s}}=\beta_{i_{s}-1}\zeta_{i_{s}-1}.
$$
Since $X_{0}^*X_{0}A=AX_{0}^*X_{0}$, $A$ leaves each Hilbert space
$\Hh_{0\ldots0,i_{s}}$ invariant.
From $X_{0}A=AX_{0}$, it follows that the restrictions of $A$ to
$\Hh_{0\ldots0,i_{s}}$ and $\Hh_{0\ldots0, i_{s}-1}$ coincide, hence
$A$ acts on $\Hh_{0\ldots0,1}$ by
$ A\zeta_{i_{s}}= A_0\zeta_{i_{s}}$, where $A_0$ is a self-adjoint operator
on $\Hh_{0\ldots0,1}$  such that
$\sigma(A_0)\sqsubseteq (q,1]$. Inserting the expression for $X_0$ into
$Z_{0\ldots0}=X_0A$ shows that
$z_{s}$ acts on $\Hh_{i_n\ldots i_{s+1}, i_{s}}$ by 
$$
z_{s}
\eta_{i_n\ldots i_{s+1},i_{s}}=q^{-(i_{s+1}+\ldots+i_{k})
+(i_{n-m+1}+\ldots+i_n)}
\beta_{i_{s}-1}A\eta_{i_n\ldots i_{s+1}, i_{s}-1}.
$$

Recall that $vZ_{0\ldots0}=Z_{0\ldots0}v$ when $l\neq 0$.
Moreover, $vZ_{0\ldots0}^*=Z_{0\ldots0}^*v$ since
$X_0v=vX_0$ and $X_0^*X_0v=vX_0^*X_0$.
Therefore it follows by the same reasoning as for $A$ that $v$
acts on $\Hh_{0\ldots0,1}$ by $v\zeta_{i_{s}}=v_0\zeta_{i_{s}}$,
where $v_0$ is a unitary operator on $\Hh_{0\ldots0,1}$.
From $vA=Av$, we conclude $v_0 A_0=A_0 v_0$.
Inserting $A_0$ and $v_0$ (when $l\neq0$) into the expressions for
$z_n,\ldots,z_{s}$, we obtain for both
$l=0$ and $l>0$ a rep\-re\-sen\-ta\-tion of the list $(m,l,k+1)$.
This completes the reduction procedure.
%

%
%XYZ
%------------------------------------------------------------------------
\subsection{Invariant integration on the quantum $(n,1)$-matrix ball }
%------------------------------------------------------------------------
                                                        \label{sec-int}
Throughout this subsection, we assume that we are given an admissible 
*-rep\-re\-sen\-ta\-tion $\pi:\qm\rightarrow\ld$ of 
the series $(m,0,k)$ such that $|Q_j|^{-1/2}\in\ld$ 
for all $j=1,\ldots, n$. 
Set $\epsilon_j=1$ if $j>k$ and $\epsilon_j=-1$ if $j\leq k$. 
Then we have $Q_j=\epsilon_j|Q_j|$. 

To develop an invariant integration theory on the 
quantum $(n,1)$-matrix ball, we proceed as in Subsection \ref{sec-qd}. 
The crucial step is to find an operator expansion of the action. 
To begin, we prove some useful operator relations. 
\begin{thl}                                         \label{Arhoz}
Define 
\begin{align}
l<n:\ \,
 \rho_l &= |Q_l|^{1/2} |Q_{l+1}|^{-1}|Q_{l+2}|^{1/2},\quad & 
\rho_n&=|Q_1|^{1/2}|Q_n|^{1/2},                   \label{defrho}\\
 A_l&=-q^{-5/2}\lambda^{-1}Q_{l+1}^{-1}z_{l+1}^*z_l, &
A_n&= q^{-1/2}\lambda^{-1}z_n,                     \label{defA}\\
 B_l &= \rho_l^{-1}A_l^*, & 
 B_n&=- \rho_n^{-1}A_n^*.                           \label{defB}
\end{align}
The operators $\rho_l$, $A_l$, and $B_l$ satisfy the following 
commutation relations:
\begin{equation}                                 \label{AB1}
\rho_i\rho_j=\rho_j\rho_i,\ \, \rho_j^{-1}\rho_j=\rho_j\rho_j^{-1}=1,
\ \,\rho_iA_j=q^{a_{ij}}A_j\rho_i,\ \,\rho_iB_j=q^{-a_{ij}}B_j\rho_i,
\end{equation}
\begin{equation}                                 \label{AB2}
A_iA_j-A_jA_i=0,\ \, i\neq j\pm1,
\quad
A_j^2A_{j\pm1}-(q+q^{-1})A_jA_{j\pm1}A_j+A_{j\pm1}A_j^2=0,
\end{equation}
\begin{equation}                                \label{AB3}
B_iB_j-B_jB_i=0,\ \,i\neq j\pm1,
\quad
B_j^2B_{j\pm1}-(q+q^{-1})B_jB_{j\pm1}B_j+B_{j\pm1}B_j^2=0,
\end{equation}
\begin{equation}                                \label{AB4}
  A_iB_j-A_jB_i=0,\ \,\, i\neq j,  \ \quad\quad
  A_jB_j-B_jA_j=\lambda^{-1}
  (\epsilon_{j+2}\epsilon_j\rho_j-\rho_j^{-1}),\,\  j<n,
\end{equation}
\begin{equation}                               \label{AB5}
         A_nB_n -B_nA_n \, = \,   -\lambda^{-1}\rho_n^{-1},
\end{equation}
where $(a_{ij})_{i,j=1}^n$ denotes the Cartan matrix of 
$sl(n+1,\C)$. 
\end{thl}
{\bf Proof.} 
Since the rep\-re\-sen\-ta\-tion is assumed to be admissible, we conclude from 
\rf[Qz] and \rf[Qz*] that  
\begin{equation}                                   \label{Q1/2z}
|Q|^{1/2}_lz_j=z_j|Q|^{1/2}_l,\,\quad|
Q|^{1/2}_lz_j^*=z_j^*|Q|^{1/2}_l,\,\quad j<l,
\end{equation}
\begin{equation}                                    \label{Q1/2z*}
|Q|^{1/2}_lz_j=qz_j|Q|^{1/2}_l,\,\quad
|Q|^{1/2}_lz_j^*=q^{-1}z_j^*|Q|^{1/2}_l,\,\quad  j\geq l.
\end{equation}
Now \rf[AB1] follows immediately from  \rf[QQ], \rf[Q1/2z], and
\rf[Q1/2z*].  The first equations of \rf[AB2]--\rf[AB4] are easily 
shown by repeated application of the commutation rules 
in $\qm$ and Equations \rf[Q1/2z] and \rf[Q1/2z*]. 
We continue with the second equation of \rf[AB4] and  \rf[AB5]. 
Using \rf[zQQ], we compute  
\begin{align*}
A_lB_l -B_lA_l &= q^{-5}\lambda^{-2}\rho_l^{-1} 
(q^2Q_{l+1}^{-1}z_{l+1}^*z_l z_l^*z_{l+1}Q_{l+1}^{-1} - 
z_l^*z_{l+1}Q_{l+1}^{-2}z_{l+1}^*z_l)                   \\
&=q^{-1}\lambda^{-2}\rho_l^{-1}Q_{l+1}^{-2}
\big( (Q_{l+2}-q^2Q_{l+1})(Q_{l+1}-Q_l)                  \\
&\quad\quad       -(Q_{l+2}-Q_{l+1})(Q_{l+1}-q^2Q_{l})\big)  \\
&=\lambda^{-1}\rho_l^{-1}(Q_{l+2}Q_{l+1}^{-2}Q_l-1)
  = \lambda^{-1}(\epsilon_{l+2}\epsilon_l\rho_l-\rho_l^{-1}),\\
A_nB_n -B_nA_n &=   
- q^{-1}\lambda^{-2}\rho_n^{-1}(q^2z_nz_n^*-z_n^*z_n)
          = -\lambda^{-1}\rho_n^{-1}.
\end{align*}
Next, we claim that
\begin{eqnarray}                                \label{AA}
A_{l-1}A_l&=&qA_lA_{l-1}-q^{-3}\lambda^{-1}Q_l^{-1}z_{l+1}^*z_{l-1},\,\ 
1<l<n,\\                                         \label{AAn}
A_{n-1}A_n&=& qA_nA_{n-1}+q^{-1}\lambda^{-1}Q_n^{-1}z_{n-1}.
\end{eqnarray}
Indeed, inserting the definition of $A_l$ and applying \rf[ball1], 
\rf[ball2], \rf[zzQ], \rf[Qz], and \rf[Qz*], one obtains 
\begin{eqnarray*}
A_{l-1}A_l&=&q^{-5}\lambda^{-2}
  Q_{l}^{-1}z_{l}^*z_{l-1}Q_{l+1}^{-1}z_{l+1}^*z_l  \\
         &=& q^{-6}\lambda^{-2}Q_{l+1}^{-1} Q_{l}^{-1}
             z_{l+1}^*(q^2z_l z_l^*-(1-q^2)Q_{l+1})z_{l-1}\\
         &=& qA_lA_{l-1}-q^{-3}\lambda^{-1}Q_l^{-1}z_{l+1}^*z_{l-1}      
\end{eqnarray*}
which proves \rf[AA]. Equation \rf[AAn] is proved similarly.
Let $l<n$. 
Multiplying \rf[AA] by $-q^{-1}A_l$ from the left and by $A_l$ from the 
right and summing both results yields \rf[AB2] with the minus sign since 
$A_lQ_l^{-1}z_{l+1}^*z_{l-1}=qQ_l^{-1}z_{l+1}^*z_{l-1}A_l$. 
Equation \rf[AB2] with the plus sign is obtained similarly by 
computing $A_{l-1}\cdot\rf[AA]-q^{-1}\rf[AA]\cdot A_{l-1}$, using 
$A_{l-1}Q_l^{-1}z_{l+1}^*z_{l-1}=q^{-1}Q_l^{-1}z_{l+1}^*z_{l-1}A_{l-1}$
and replacing $l$ by $l+1$ for $l+1<n$.
The same steps applied to \rf[AAn] yield \rf[AB2] for $l=n$ and
$l+1=n$  since also 
$A_nQ_n^{-1}z_{n-1}=qQ_n^{-1}z_{n-1}A_n$ and 
$A_{n-1}Q_n^{-1}z_{n-1}=q^{-1}Q_n^{-1}z_{n-1}A_{n-1}$. 
The second equations of \rf[AB3] follow from the second 
equations of \rf[AB2] by applying the involution and multiplying 
by $\rho_j^{-2}\rho_{j\pm1}^{-1}$.                   \hfill $\Box$
\mn\\
{\bf Remark.} By \rf[AB5], the operators $A_l$, $B_l$, and $\rho_l$ 
do not satisfy the defining relations of $\sutn$. 
If $n>1$, then we get only for the series $(n,0,0)$ 
a *-rep\-re\-sen\-ta\-tion of $\sutm$ by assigning 
$K_j$ to $\rho_j$, $E_j$ to $A_j$, and $F_j$ to $B_j$, $j<n$. 

To see this, observe that we must have 
$\epsilon_{j+2}=\epsilon_{j}$ by \rf[AB4]. But  
$\epsilon_{n-1}=1$ since $Q_{n+1}=1$, and
$\epsilon_{n}=1$ since 
$\epsilon_n|Q_n|=Q_{n-1}+z_{n-1}z_{n-1}^*>0$
by \rf[zQQ]
(cf.\ the remarks after Proposition \ref{DarMat}), 
so $\epsilon_{n}= \ldots=\epsilon_{1}=1$.  \mn  

Although Equations \rf[AB1]--\rf[AB5] do not yield a rep\-re\-sen\-ta\-tion
of $\sutn$, the analogy to \rf[sutm1]--\rf[sutm4] is obvious, 
so it is natural to try to define an operator expansion of the 
action by imitating the adjoint action. That this can be done is 
the assertion of the next lemma. 
Again, we write $f$ instead of $\pi(f)$ and 
$X\ang f$ instead of $\pi(X\ang f)$ for $f\in\qm$, $X\in\sutm$.

\begin{thl}                               \label{ballact}
With the  operators $\rho_l$, $A_l$, and $B_l$ defined
in Lemma \ref{Arhoz}, set 
\begin{eqnarray}
  K_j\,\ang\, f &=&\rho_jf\rho_j^{-1},\ \quad K_j^{-1}\,\ang\, f
 \ =\ \rho_j^{-1}f\rho_j,                           \label{ballact1}\\
  E_j\,\ang\, f &=& A_jf-\rho_jf\rho_j^{-1}A_j,     \label{ballact2}\\
  F_j\,\ang\, f &=& B_jf\rho_j-q^2f\rho_jB_j         \label{ballact3}     
\end{eqnarray}       
for $j=1,\ldots,n$. Then Equations \rf[ballact1]--\rf[ballact3] applied to 
$f\in\qm$ define an operator expansion of the action $\ang$ on $\qm$. 
The same formulas applied to $f\in \ld$ turn the 
$O^*$-algebra $\ld$ into a $\sutn$-module *-algebra.
\end{thl}
{\bf Proof.} The lemma is proved by direct verifications. We start 
by showing that $\ld$ with the $\sutn$-action defined by 
\rf[ballact1]--\rf[ballact3] becomes a $\sutn$-module *-al\-ge\-bra. 
That the action satisfies \rf[modalg]--\rf[modstar] is readily seen  
if we replace in the proof of Lemma \ref{discact} 
$y^{\pm1}$ by $\rho_j^{\pm1}$, $A$ by $A_j$, and $B$ by $B_j$.
By using Lemma \ref{Arhoz}, it is easy to check that the action 
is consistent with \rf[sutm1] and the first relations 
of \rf[sutm2]--\rf[sutm4]. For example, \rf[AB1] and \rf[AB2] give 
\begin{align*}
(K_iE_j)\ang f&= \rho_i(A_jf-\rho_jf\rho_j^{-1}A_j)\rho_i^{-1}
         =q^{a_{ij}}(A_j\rho_i f\rho_i^{-1}
            -\rho_j\rho_i f\rho_i^{-1}\rho_j^{-1}A_j)\\
          &=(q^{a_{ij}}E_jK_i)\ang f,\\
(E_lE_j)\ang f&= A_l(A_jf-\rho_jf\rho_j^{-1}A_j) 
-\rho_l(A_jf-\rho_jf\rho_j^{-1}A_j) \rho_l^{-1}A_l\\
  &=  A_j(A_lf-\rho_lf\rho_l^{-1}A_l)
-\rho_j(A_lf-\rho_lf\rho_l^{-1}A_l)\rho_j^{-1}A_j
 =(E_lE_j)\ang f
\end{align*}
for all $f\in \ld$, $i,j,l=1,\ldots,n$, $l\neq j\pm1$. 
As in the proof of Lemma \ref{discact}, we have 
\begin{align*}
(E_jF_j-F_jE_j)\ang f &=
A_j B_jf\rho_j+\rho_jfB_jA_j-B_jA_jf\rho_j-\rho_jfA_j B_j\\
&=(A_jB_j-B_jA_j)f\rho_j - \rho_jf(A_jB_j-B_jA_j).
\end{align*}
Inserting \rf[AB4] if $j<n$ and \rf[AB5] if $j=n$  shows that 
the action is consistent with the second equation of \rf[sutm4].

We continue with the second equation of \rf[sutm2]. A straightforward 
calculation shows that 
\begin{align}                    \label{EEE1}
(E_j^2E_{j\pm1})\ang f &=
A_j^2\Aj f - \rho_j^2\rj f \rho_j^{-2}\rj^{-1}\Aj A_j^2  \\
&\phantom{=} 
- (q+q^{-1})A_j\Aj \rho_j f\rho_j^{-1}A_j   
+\Aj \rho_j^2f\rho_j^{-2}A_j^2                    \notag  \\ 
&\phantom{=}   - A_j^2\rj f \rj^{-1}\Aj 
 +(q+q^{-1}) A_j \rho_j\rj f\rho_j^{-1}\rj^{-1}\Aj A_j, \notag
\end{align}
\begin{align}                                  \label{EEE2}
(E_jE_{j\pm1}E_j)\ang f &=     
A_j\Aj A_j f - \rho_j^2\rj f \rho_j^{-2}\rj^{-1}A_j\Aj A_j\\
&\phantom{=}  - q\Aj A_j \rho_j f\rho_j^{-1}A_j
-A_j\Aj \rho_j f\rho_j^{-1}A_j        \notag     \\
&\phantom{=}  +q\Aj\rho_j^2 f\rho_j^{-2}A_j^2 
       -q^{-1}A_j^2\rj f\rj^{-1}\Aj      \notag   \\
&\phantom{=}  +A_j \rho_j\rj f\rho_j^{-1}\rj^{-1}\Aj A_j   \notag
\hspace{-0.1pt}+\hspace{-0.1pt}
q^{-1}A_j \rho_j\rj f\rho_j^{-1}\rj^{-1}A_j\Aj,
\end{align}
\begin{align}                   \label{EEE3}
(E_{j\pm1}E_j^2)\ang f &=          
\Aj A_j^2 f - \rho_j^2\rj f \rho_j^{-2}\rj^{-1}A_j^2\Aj \\
&\phantom{=} -(1+q^{-2})\Aj A_j \rho_j f\rho_j^{-1}A_j     
+q^{2}\Aj\rho_j^2 f\rho_j^{-2}A_j^2           \notag  \\ 
&\phantom{=} -q^{-2}A_j^2\rj f\rj^{-1}\Aj 
 +(1+q^{-2})A_j \rho_j\rj f\rho_j^{-1}\rj^{-1}A_j\Aj, \notag
\end{align}
where we repeatedly used \rf[AB1].
Taking the sum $\rf[EEE1]-(q+q^{-1})\cdot\rf[EEE2]+\rf[EEE3]$ gives $0$ 
since the sums over the first and the second summands vanish 
by \rf[AB2] and the other summands cancel. 
The last result implies also the second relation of \rf[sutm3] 
since $X\ang f= (S(X)^*\ang f)^*$ for all $X\in\sutn$ 
and $S(F_j)^*=-(-1)^{\delta_{nj}}q^2E_j$.    

It remains to prove that \rf[ballact1]--\rf[ballact3] define an operator 
expansion of the action. 
That \rf[ballact1] yields the action of $K^{\pm1}$ on $z$ and $z^*$ 
is easily verified by using \rf[Q1/2z] and \rf[Q1/2z*]. 
Let $l<n$. 
Applying \rf[ball1], \rf[ball2], \rf[Q1/2z], and \rf[Q1/2z*], we get 
$\rho_lz_j\rho_l^{-1}A_l=A_lz_j$ and 
$\rho_lz_j^*\rho_l^{-1}A_l=A_lz^*_j$ whenever $j\notin\{l,l+1\}$, 
hence $E_l\ang z_j=E_l\ang z_j^*=0$. 
Similarly, $\rho_lz_l\rho_l^{-1}A_l=qz_lA_l=A_lz_l$ and 
$\rho_lz_{l+1}^*\rho_l^{-1}A_l=q^{-1}z_{l+1}^*A_l=A_lz_{l+1}^*$, 
so $E_l\ang z_{l}=E_l\ang z_{l+1}^*=0$. 
Equation \rf[ballact2] applied to $z_{l+1}$ and $z_l^*$ gives 
\begin{align*}
E_l\ang z_{l+1}&= A_lz_{l+1}-q^{-1}z_{l+1}A_l
=-q^{-3/2}\lambda^{-1}Q_{l+1}^{-1}(z_{l+1}^*z_{l+1}-z_{l+1}z_{l+1}^*)z_l\\
&=q^{-1/2}z_l,
\end{align*}
$$
E_l\ang z_{l}^* = A_lz_{l}^*-q^{-1}z_{l}^*A_l
=-q^{-5/2}\lambda^{-1}z_{l+1}^*
(q^2 z_{l}z_{l}^*-z_{l}^*z_{l})Q_{l+1}^{-1}
=-q^{3/2}z_{l+1}^*,
$$
where we used \rf[zzQ]. For $E_n$, we obtain 
\begin{align*}
E_n\ang z_j&=q^{-1/2}\lambda^{-1}(z_nz_j-qz_jz_n)= -q^{1/2}z_nz_j,\\
E_n\ang z_j^*&=q^{-1/2}\lambda^{-1}(z_nz_j^*-q^{-1}z_j^*z_n)=0
\end{align*}
if $j<n$. The action of $E_n$ on $z_n$ and $z_n^*$ is calculated 
analogously to \rf[Ez] and \rf[Ez*]. We have thus proved that 
\rf[ballact2] is consistent with the action of $E_j$, $j=1,\ldots,n$,  
on $\qm$. The corresponding result for $F_j$ follows from this by 
using $F_j\ang f=-(-1)^{\delta_{nj}}q^2(E_j\ang f^*)^*$. 
                                                      \hfill $\Box$
\mn

Let $\omega_1,\ldots,\omega_n$ be the simple roots of the 
Lie algebra $sl_{n+1}$. For $\gamma = \sum_{j=1}^n p_j\omega_j$,  
we write $K_\gamma=K_1^{p_1}\cdots K_n^{p_n}$. Recall that, for a 
finite dimensional rep\-re\-sen\-ta\-tion $\sigma$ 
of $\sutn$, the quantum trace 
$$
           \tr[q,L]a:= \tr\sigma(aK_{2\omega}^{-1})
$$ 
defines an invariant integral on $\sutn$, where $\omega$ denotes 
the half-sum of all positive roots (see \cite[Proposition 7.14]{KS}). 
$K_{2\omega}$ is chosen such that 
$X K_{2\omega}=K_{2\omega}S^2(X)$ for all  $X\in\sutn$. 
In Subsection \ref{sec-qd}, we replaced $K\,(=K_{2\omega})$ by $y$ and 
proved the existence of  invariant integrals on 
appropriate classes of functions. Our aim is to generalize this result 
to $\qm$. 

The half-sum of positive roots is given by  
$\omega=\frac{1}{2}\sum_{l=1}^nl(n-l+1)\omega_l$. 
Consider $\Gamma:=\prod_{l=1}^n \rho_l^{-l(n-l+1)}$. 
Inserting the definition of $\rho_l$ gives  
\begin{equation}                                   \label{Gamma}
 \Gamma=|Q_1|^{-n}|Q_2|\cdots |Q_n|,\,\ n>1,\, \quad 
 \Gamma=|Q_1|^{-1},\,\ n=1,
\end{equation}
since $-\frac{1}{2}(l-1)(n-l+2)+l(n-l+1)-\frac{1}{2}(l+1)(n-l)=1$
for $1<l\leq n$. The operator $|Q_1|$ appears
in the definition of $\Gamma$ twice, in $\rho_1^{-n}$ 
and $\rho_n^{-n}$, in each factor to the power $-n/2$. 
For $n=1$, Equation \rf[Gamma] is trivial (cf.\ Equation \rf[defrho]). 
The following proposition shows that $\Gamma$ enables us  
to define a generalization of the quantum trace. 

Notice that $z_n$, $z_n^*$, $K_n^{\pm1}$, $E_n$, and $F_n$ satisfy the 
relations of the quantum disc, in particular, Equation \rf[h=0] 
applies. Therefore we cannot have a normalized invariant integral 
on $\qm$. 
\begin{thp}                                       \label{haarball}
Let $\fA$ be the O*-algebra generated by the operators 
$z_j$,\ $z_j^*$,\ $|Q_j|^{1/2}$, and $|Q_j|^{-1/2}$, 
$j=1,\ldots,n$. 
Then the *-algebras $\FFD$ and $\BBA$ defined in 
\rf[FD] and \rf[BA], respectively, are $\sutn$-module *-algebras,
where the action is given by \rf[ballact1]--\rf[ballact3]. 
The linear functional
\begin{equation}                                  \label{bh}
    h(f):=c\, \tr \ov{f\Gamma},\quad c\in\R,        
\end{equation}
defines an invariant integral on both $\FFD$ and $\BBA$. 
\end{thp}
{\bf Proof.} From the definition of $\FFD$ and $\BBA$, it is 
obvious that both algebras are stable under the $\sutn$-action
defined by \rf[ballact1]--\rf[ballact3], in particular, by 
Lemma \ref{ballact}, they are $\sutn$-module *-algebras. 

We proceed as in the proof of Proposition \ref{haar} and show 
the invariance of $h$ for generators by using the trace property 
$\tr \ov{agb}=\tr \ov{gba}=\tr \ov{bag}$ for all 
$g\in\BBA$, $a,b\in\fA$. Let $g\in\BBA$. 
Clearly, $\rho_l$ commutes with $\Gamma$, hence 
$$
h(K_l^{\pm1}\ang g)=\tr\ov{\rho_l^{\pm1}g\rho_l^{\mp1}\Gamma}=
 \tr\ov{g\Gamma}=\vare(K_l^{\pm1})h(g)
$$
It follows from the definition of $\Gamma$ and from \rf[AB1] 
that $A_l\Gamma=q^2\Gamma A_l$ for all $l$ since 
$$
-(l-1)(n-l+2)+2l(n-l+1)-(l+1)(n-l)=2.
$$ 
Hence  
$\rho_l^{-1}A_l\Gamma=\Gamma A_l\rho_l^{-1}$ 
and therefore 
$$
h(E_l\ang g)= 
\tr(\ov{A_lg\Gamma}-\ov{\rho_lg\rho_l^{-1}A_l\Gamma})=
\tr\ov{A_lg\Gamma}-\tr\ov{A_lg\Gamma}=0=\vare(E_l)h(g). 
$$
Applying the involution to $A_l\Gamma=q^2\Gamma A_l$ 
shows that $\Gamma B_l=q^2B_l\Gamma$. Thus 
$$
h(F_l\ang g)= 
\tr(\ov{B_lg\rho_l\Gamma}-q^2\ov{g\rho_lB_l\Gamma})=
\tr\ov{B_lg\rho_l\Gamma}-\tr\ov{B_l\rho_lg\Gamma}=0
  =\vare(F_l)h(g). \\
$$
This completes the proof.                     \hfill $\Box$
\mn\\
{\bf Remark.}
As in Subsection \ref{sec-qd}, we consider $\BBA$ as 
the algebra of infinitely differentiable functions which vanish 
sufficiently rapidly at ``infinity'' and $\FFD$ as the 
infinitely differentiable functions with compact support. 
%
%
%
%------------------------------------------------------------------------
\subsection{Topological aspects of *-representations}
%------------------------------------------------------------------------
%
In this subsection, we shall restrict ourselves to 
irreducible *-rep\-re\-sen\-ta\-tions of the series $(m,0,k)$. 
Let $D$ denote the linear space defined in Proposition \ref{DarMat}. 
Then the operators $|Q_j|^{-1/2}$, $1\le j\le n$, belong to 
$\ld$ and the O*-algebra $\fA$ of Proposition \ref{haarball} is 
well defined. As in Subsection \ref{topasp}, we prefer for topological 
reasons to work with closed O*-algebras. In particular, 
we suppose that the *-rep\-re\-sen\-ta\-tion is given on the domain 
$D_{\fA}=\cap_{a\in\fA}D(\bar a)$. It turns out that the topological 
properties of the (closed) O*-algebra $\fA$ are very similar to that of 
Subsection \ref{topasp}. 
\begin{thl} \mbox{ }                                 \label{l-topball}
\begin{enumerate}
\item
  $\fA$ is a commutatively dominated O*-algebra on a 
  Frechet domain.
\item 
   $D_\fA$ is nuclear, in particular,
   $D_\fA$ is a Frechet--Montel space.
\end{enumerate}
\end{thl}
{\bf Proof.}
The operator 
\begin{equation}
  T:=1+Q_1^2+\ldots +Q_n^2+Q_1^{-2}+\ldots+Q_n^{-2}    \label{T}
\end{equation}
is essentially self-adjoint on $D_\fA$ , and $T>2$. Let $\varphi\in D_\fA$.
As in the proof of 
Lemma \ref{l-top}, we conclude from a standard argument that, for 
each polynomial 
$p=p(|Q_1|^{1/2},\ldots,|Q_n|^{1/2},|Q_1|^{-1/2},\ldots,|Q_n|^{-1/2})$,  
there exist  $k\in \N$ such that 
$||p\varphi||\le||T^k\varphi||$. Furthermore, for each finite sequence
$k_1,\ldots,k_N\in\N$ 
and real numbers $\gamma_1,\ldots, \gamma_N\in(0,\infty)$, 
we find $k_0\in\N$ such that 
$\sum_{j=1}^N\gamma_j ||T^{k_j}\varphi||\le ||T^{k_0}\varphi||$. 
Let $p$ be as above and let 
$I,J\in\N^n$ such that $I\cdot J=0$. By \rf[ball1]--\rf[ball4] and 
\rf[zQQ]--\rf[QQ], $(z^Ipz^{*J})^*(z^Ipz^{*J})$ is a polynomial 
in $|Q_j|^{1/2}$, $|Q_j|^{-1/2}$, $j=1,\ldots,n$, say $\tilde p$. 
Thus  there exist $k\in\N$ such that 
$$
||z^Ipz^{*J}\varphi||=\ip{\tilde p\varphi}{\varphi}^{1/2}
\le(||\tilde p\varphi||\,||\varphi||)^{1/2}\le||T^k\varphi||.   
$$
From the definition of $\fA$, \rf[f], \rf[Q1/2z], and \rf[Q1/2z*],  
it follows that each $f\in\fA$ can be written 
as $f=\sum_{I\cdot J=0}z^Ip_{IJ}z^{*J}$, where $p_{IJ}$ are polynomials 
in $|Q_j|^{1/2}$, $|Q_j|^{-1/2}$, $j=1,\ldots,n$. 
From the preceding arguments, we conclude  
that there exist  $m\in\N$ such that 
$||f\varphi||\le ||T^m\varphi||$ for all $\varphi\in D_\fA$, 
therefore $||\cdot||_f\le||\cdot||_{T^m}$.  
This implies that the family 
$\{ ||\cdot||_{T^{2^k}}\}_{k\in\N}$ generates the graph topology 
on $D_\fA$ and 
$D_\fA=\cap_{k\in\N}D(\bar{T}^{2^k})$ which proves (i).

Note that the proof of Lemma \ref{l-top}(ii) is based on the 
observation that the operator $\bar T^{-1}$ is a Hilbert--Schmidt 
operator. One 
easily checks that this holds also for the operator $T$ defined
in (\ref{T}). Now the rest of the proof runs completely 
analogous to that of Lemma \ref{l-top}.
                                          \hfill $\Box$
\begin{thp}  \mbox{ }                  \label{topball}
\begin{enumerate}
\item
   $\FFDA$ is dense in $\ldda$ with respect to the bounded 
   topology $\tau_b$.
\item
   The $\sutn$-action on $\ldop$ is continuous with respect to $\tau_b$.  
\end{enumerate}
\end{thp}

The proof of Proposition \ref{topball} is completely analogous
to that of Proposition \ref{top}. 

\begin{thc}
Let $\FFD$ denote the O*-algebra of finite rank operators
on $D$ defined in \rf[FD]. 
Then $\FFD$ is a $\sutn$-module *-subalgebra of $\FFDA$ 
and $\FFD$ is dense in $\ldda$. 
\end{thc}
{\bf Proof.} 
Since $D\subset D_\fA$, we can consider $\FFD$ as a *-subalgebra 
of $\FFDA$. It follows from Proposition \ref{DarMat} that 
$\FFD$ is stable under the $\sutn$-action defined in 
Lemma \ref{ballact}, in particular, it is a 
$\sutn$-module *-algebra. The density of $\FFD$ in $\ldda$ can be proved 
in exactly the same way as in  Corollary \ref{dense}. 
                                           \hfill $\Box$
\mn
  
Recall that the self-adjoint operators $\bar Q_j$, $j=1,\ldots,n$, strongly 
commute. Set 
$$
  \gM:=\sigma(\bar Q_1)\times\ldots\times
        \sigma(\bar Q_n).
$$
By the spectral theorem of self-adjoint operators, 
we can assign to each (Borel measurable) function 
$\psi : \gM \rightarrow \C$ an operator 
$\psi(\bar Q_1,\ldots,\bar Q_n)$ such that 
%%%%%%%%%for type $(m,0,k)$ rep\-re\-sen\-ta\-tions 
$$
\psi(\bar Q_1,\ldots,\bar Q_n)\eta_{i_n\ldots i_1}=
\psi(t_{i_1},\ldots,t_{i_n})\eta_{i_n\ldots i_1},
$$ 
where $t_{i_j}=q^{2(i_j+\ldots+i_n+\alpha)}$ for $j>k$,
$t_{i_j}=-q^{-2(i_j+\ldots+i_k)+2(i_{k+1}+\ldots+i_n+\alpha)}$ 
for $j\le k$, and $A=q^{2\alpha}$. ($A$ denotes the operator 
appearing in the type $(m,0,k)$ rep\-re\-sen\-ta\-tions for $k>0$.
If $k=0$, set $\alpha=0$.)\, 
Define 
\begin{align*}
\lefteqn{\cS(\gM)=}\\
&\quad \{\,\psi:\gM\rightarrow\C\,;\sup_{(t_1,\ldots,t_n)\in\gM}
|t_1^{s_1}\cdots t_n^{s_n}\psi(t_1,\ldots,t_n)|<\infty\,\ 
\mathrm{for\ all}\ s_1,\ldots,s_n\in\Z\,\} 
\end{align*}
and
$$
\cS(D)=\{\,\sum_{I\cdot J=0} z^I\psi_{IJ}(\bar Q_1,\ldots,\bar Q_n)z^{*J}
\,;\, \psi_{IJ}\in \cS(\gM),\ \#\{\psi_{IJ}\neq0\}<\infty\,\}. 
$$
\begin{thl}                         \label{SM}
With the action defined in Lemma \ref{ballact}, $\cS(D)$ becomes a 
$\sutn$-module *-subalgebra of $\BBA$.  
The operators $z_j$, $z_j^*$, $j=1,\ldots,n$,  
and $\psi(\bar Q_1,\ldots,\bar Q_n)$, $\psi\in\cS(\gM)$, 
satisfy the following commutation rules
$$
  \psi(\bar Q_1,\ldots,\bar Q_j,\bar Q_{j+1},\ldots,\bar Q_n)z_j=
   z_j\psi(q^2\bar Q_1,\ldots,q^2\bar Q_j,\bar Q_{j+1},\ldots,\bar Q_n)
$$
$$
 z_j^* \psi(\bar Q_1,\ldots,\bar Q_j,\bar Q_{j+1},\ldots,\bar Q_n)=
 \psi(q^2\bar Q_1,\ldots,q^2\bar Q_j,\bar Q_{j+1},\ldots,\bar Q_n)z_j^*.
$$
\end{thl}

The proof of Lemma \ref{SM} differs from that of Lemma \ref{SD} 
only in notation, the argumentation to establish the result 
remains the same. 

Since $\FFD\subset\FFDA$ and $\cS(D)\subset\BBA$, we can consider 
$\FFD$ and $\BBA$ as algebras of  infinite differentiable 
functions with compact support and 
which are rapidly decreasing, respectively. 
It is not difficult to see that $\FFD$ is the set of all
$\sum_{I\cdot J=0}z^I\psi_{IJ}(\bar Q_1,\ldots,\bar Q_n)z^{*J}\in\cS(D)$, 
where the functions $\psi_{IJ}\in\cS(\gM)$ have 
finite support. On $\cS(D)$, we have the following explicit formula of 
the invariant integral. 
\begin{thp}                                        \label{hsutn}
Set
$$
  \gM_0:=\sigma(\bar Q_1)\backslash\{0\}\times\ldots\times
        \sigma(\bar Q_n)\backslash\{0\}.
$$
Assume that 
$f=\sum_{I\cdot J=0}z^I\psi_{IJ}(\bar Q_1,\ldots,\bar Q_n)z^{*J}\in\cS(D)$.
Then the invariant integral $h$ defined in Proposition \ref{haarball} 
is given by
$$
 h(f)= c \sum_{(t_1,\ldots,t_n)\in\gM_0}\psi_{00}(t_1,\ldots,t_n)
|t_1|^{-n}|t_2|\cdots|t_n|. 
$$
(If $n=1$, then $t_2,\ldots,t_n$ are omitted.) 
\end{thp}
{\bf Proof.} Recall that $h(f)=c\, \tr \ov{f\Gamma }$, where $\Gamma$ 
is given by \rf[Gamma]. If $I\neq(0,\ldots,0)$ or $J\neq(0,\ldots,0)$, 
then 
$\ip{\eta_{i_n\ldots i_1}}{z^I\psi_{IJ}(\bar Q_1,\ldots,\bar Q_n)z^{*J}
\Gamma \eta_{i_n\ldots i_1}}=0$
since, by Proposition \ref{DarMat}, $\psi_{IJ}(\bar Q_1,\ldots,\bar Q_n)$ 
and $\Gamma$ are diagonal and $z^I$ and $z^{*J}$ act as shift operator 
on $\Hh$. 
Hence only 
$\psi_{00}(\bar Q_1,\ldots,\bar Q_n)\Gamma $ contributes to the 
trace. 

For each tuple $(t_1,\ldots,t_n)\in \gM_0$, there exists exactly 
one tuple $(i_1,\ldots,i_n)$ such that $\eta_{i_n\ldots i_1}\in D$ and
$Q_j \eta_{i_n\ldots i_1}=t_j\eta_{i_n\ldots i_1}$, $j=1,\ldots,n$. 
This can be seen inductively; $Q_n$ determines $i_n$ uniquely, and 
if $i_n,\ldots,i_{n-k+1}$ are fixed, then $Q_{n-k}$ determines 
uniquely $i_{n-k}$ (see the remark after Proposition \ref{DarMat}). 
Since the vectors $\eta_{i_n\ldots i_1}$ constitute an orthonormal basis 
of eigenvectors of the $Q_j$'s, 
and since $\Gamma$ is given by  
$\Gamma=|Q_1|^{n}|Q_2|^{-1}\cdots |Q_n|^{-1}$ for $n>1$, 
$\Gamma=|Q_1|$ for $n=1$, the assertion follows. 
                                                   \hfill $\Box$
\mn

In the following, let  $n>1$. 
We noted in Subsection \ref{sec-ball} that the action of the 
elements $E_j$, $F_j$, $K_j^{\pm1}$, $j=1,\ldots,n-1$ on $\qm$ 
induces  a $\sutm$-action 
which turns $\qm$ into a $\sutm$-module *-algebra. 
$\sutm$ is regarded as a compact real form of 
$\cU_q(\mathrm{sl}_{n})$. Naturally, the compactness should 
be manifested in the existence of a normalized invariant 
integral on $\qm$. This is indeed the case. Consider a 
irreducible *-rep\-re\-sen\-ta\-tion of type $(n,0,0)$. Then the 
operators $Q_j$, $j=1,\ldots,n$, are bounded, and $Q_1$ 
is of trace class. In Proposition \ref{haarball}, a 
$\sutn$-invariant functional $h$ was given by 
$h(f):=c\, \tr \ov{f\Gamma }$, where 
$\Gamma =|Q_1|^{-n}|Q_2|\cdots |Q_n|$. 
Notice that the proof of Proposition \ref{haarball} uses 
only the commutation relations of $\Gamma $ with 
$A_i$, $B_i$, and $\rho_i$, $i=1,\ldots,n$. The crucial observation 
is that $Q_1$ commutes with 
$A_j$, $B_j$, and $\rho_j$, $j=1,\ldots,n-1$. Therefore 
the commutation relations used in proving the invariance of $h$ 
remain unchanged if we multiply $\Gamma $ by 
$Q_1^{n+1}$. Furthermore, $\Gamma Q_1^{n+1}$ is of 
trace class. This suggests that 
$h(f):=c\, \tr f\Gamma Q_1^{n+1}$ defines a $\sutm$-invariant 
integral on $\qm$. The only difficulty is that the definitions of 
$A_j$, $B_j$, and $\rho_j^{\pm1}$ involve the unbounded operators 
$Q_j^{-1}$, therefore we cannot freely apply the trace property 
in proving the invariance of $h$. Nevertheless, 
a modified proof will establish the result. 
\begin{thp}                                  \label{hqm}
Let $n>1$ and set $c:=\prod_{k=1}^n(1-q^{2k})^{-1}$. 
Suppose we are given an irreducible *-rep\-re\-sen\-ta\-tion 
of $\qm$ of type $(n,0,0)$. 
Then the linear functional 
\begin{equation}                                  \label{bbh}
    h(f):=c\, \tr f\Gamma Q_1^{n+1}
     =c\, \tr fQ_1\cdots Q_n\,,\,\quad f\in\qm,        
\end{equation}
defines a normalized\, $\sutm$-invariant integral on $\qm$.  
\end{thp}

{\bf Proof.} First note that the vectors $\eta_{i_n\ldots i_1}$, 
$i_1,\ldots,i_n\in\N_0$, form a complete set of eigenvectors of the 
positive operator $Q_1$ with corresponding eigenvalues 
$q^{2(i_1+\ldots+i_n)}$. 
As $\sum_{i_1,\ldots,i_n\in\N_0}q^{2(i_1+\ldots+i_n)}<\infty$, 
$Q_1$ is of trace class. This implies that 
$f\Gamma Q_1^{n+1}=fQ_1\cdots Q_n$ is of trace class for all 
$f\in\qm$ since the rep\-re\-sen\-ta\-tions of the series $(n,0,0)$ are 
bounded. Therefore $h$ is well defined. 
An easy calculation  shows that $h(1)=1$. 

As in the proof of Proposition \ref{haar}, it suffices to verify the  
invariance of $h$ for the generators of $\sutm$.  
Recall that $\qm$ is the linear span of the elements 
$z^Ip_{IJ}z^{*J}$, where ${I,J\in\N_0^n}$, $I\cdot J=0$, and
$p_{IJ}$ is a polynomial in $Q_i$, $i=1,\ldots,n$. 
If $I\neq0$ or $J\neq0$, then
\begin{equation}                                  \label{K0}
0=\ip{\eta_{i_n\ldots i_1}}{\rho_j^{\pm1}z^Ip_{IJ}z^{*J}\rho_j^{\mp1}
\Gamma Q_1^{n+1}\eta_{i_n\ldots i_1}}=
\ip{\eta_{i_n\ldots i_1}}{z^Ip_{IJ}z^{*J}
\Gamma Q_1^{n+1}\eta_{i_n\ldots i_1}}
\end{equation}
by the same arguments as in the proof of Proposition \ref{hsutn}.
Hence 
$$
h(K_j^{\pm1}\ang(z^Ip_{IJ}z^{*J}))=
\vare(K_j^{\pm1})h(z^Ip_{IJ}z^{*J})=0.
$$ 
If $I=J=0$, then 
$$
K_j^{\pm1}\ang p_{IJ}=\rho_j^{\pm1}p_{IJ}\rho_j^{\mp1}=p_{IJ},
$$ 
thus $h(K_j^{\pm1}\ang p_{IJ})=h(p_{IJ})=\vare(K_j^{\pm1})h(p_{IJ})$. 
This proves the invariance of $h$ with respect to $K_j^{\pm1}$, 
$j=1,\ldots, n-1$.  

Recall that $A_k=-q^{-5/2}\lambda^{-1}Q_{k+1}^{-1}z_{k+1}^*z_k$, $k<n$.  
If $I\neq(0,\ldots,1,\ldots,0)$ or $J\neq(0,\ldots,1,\ldots,0)$
with $1$ in the $(k+1)$th and $k$th positions, respectively, then
we have similarly to Equation \rf[K0] 
\begin{eqnarray*}
0 &=& \ip{\eta_{i_n\ldots i_1}}{A_kz^Ip_{IJ}z^{*J}
\Gamma Q_1^{n+1}\eta_{i_n\ldots i_1}}   \\
&=& \ip{\eta_{i_n\ldots i_1}}{\rho_kz^Ip_{IJ}z^{*J}\rho_k^{-1}A_k
\Gamma Q_1^{n+1}\eta_{i_n\ldots i_1}}. 
\end{eqnarray*}
Thus $h(E_k\ang(z^Ip_{IJ}z^{*J}))=
\vare(E_k)h(z^Ip_{IJ}z^{*J})=0$. 

Now let $p$ denote an arbitrary polynomial in $Q_i$, $i=1,\ldots,n$. 
Then, by the definition of $\Gamma$ and  repeated application of the 
commutation rules of $Q_i$ 
with $z_j$ and  $z_j^*$, we obtain 
\begin{align*}
Q_{k+1}^{-1}z_{k+1}^*z_kz_{k+1}pz_k^*\Gamma Q_1^{n+1}&=
z_{k+1}^*z_kz_{k+1}pz_k^*Q_1\cdots Q_kQ_{k+2}\cdots Q_n,\\ 
\rho_kz_{k+1}pz_k^*\rho_k^{-1}Q_{k+1}^{-1}z_{k+1}^*z_k\Gamma 
Q_1^{n+1}&=z_{k+1}pz_k^*Q_1\cdots Q_kQ_{k+2}\cdots Q_nz_{k+1}^*z_k.
\end{align*}
All operators on the right hand sides are bounded and $Q_1$ is of 
trace class, in particular, the trace property applies.   
Therefore, the difference of the traces of the right hand sides 
vanishes. Hence 
\begin{align*}
h(E_k\ang (z_{k+1}pz_k^*))&=c\,\tr (A_kz_{k+1}pz_k^*-
\rho_kz_{k+1}pz_k^*\rho_k^{-1}A_k)\Gamma Q_1^{n+1}=0\\
&=\vare(E_k)h(z_{k+1}pz_k^*) 
\end{align*}
which establishes  the invariance of $h$ with respect to 
$E_k$, $k=1,\ldots,n-1$. 

To verify that $h$ is invariant  
with respect to $F_k$, $k=1,\ldots,n-1$, notice that 
$h(f^*)=\ov{h(f)}$ for all $f\in\qm$ since the operator 
$\Gamma Q_1^{n+1}$ is self-adjoint. 
Thus, by \rf[modstar] and the preceding,  
$$
h(F_k\ang f)=\ov{h(S(F_k)^{*}\ang f^*)}=-q^2\ov{h(E_k\ang f^*)}
=0=\vare(F_k)h(f) 
$$
for all $f\in\qm$.                       \hfill $\Box$
\begin{thc}                                         \label{hhqm}
Let $f=\sum_{I\cdot J=0}z^Ip_{IJ}(Q_1,\ldots,Q_n)z^{*J}\in\qm$. Then the 
invariant integral $h$ defined in Proposition \ref{hqm} is given by 
$$
 h(f)= c \sum_{j_1,\ldots,j_n\in\N_0}p_{00}(q^{j_1},\ldots,q^{j_n})
 q^{j_1}\cdots q^{j_n}. 
$$
\end{thc}
{\bf Proof.} 
Taking into account that 
$\gM_0=\{(q^{j_1},\ldots,q^{j_n})\,;\,j_1,\ldots,j_n\in\N_0\}$
for rep\-re\-sen\-ta\-tions of the series $(n,0,0)$, 
Corollary \ref{hhqm} is verified by an obvious modification of 
the proof of Proposition \ref{hsutn}.          \hfill $\Box$
%
%
%*********************************************************************
%
\section{Concluding remarks}
%
%**********************************************************************
% 
%
In general, the definition of quantum groups and quantum spaces 
is completely algebraic. 
However, our definition of integrable functions involves 
operator algebras.  
The discussion in this paper shows that 
operator algebras form a natural setting for the study 
of non-compact quantum spaces. 
For example, Hilbert space 
rep\-re\-sen\-ta\-tions provide us with the powerful 
tool of spectral theory which allows to define functions of self-adjoint 
operators. 
We emphasize that different 
rep\-re\-sen\-ta\-tions will lead to different algebras 
of integrable functions. 
If one accepts that rep\-re\-sen\-ta\-tions carry 
information about the underlying quantum space (for instance, 
by considering the spectrum of self-adjoint operators), then 
rep\-re\-sen\-ta\-tions can be used to distinguish between $q$-deformed 
manifolds which are isomorphic on purely algebraic level.

The crucial step of our approach 
was to find an operator expansion of the action. 
At first sight it seems  a serious drawback that no direct method 
was given to obtain an operator expansion of the action. 
This problem can be removed by considering cross product algebras. 
Inside the cross product algebra, the action can be expressed by 
algebraic relations. Rep\-re\-sen\-ta\-tions of cross product algebras 
lead therefore to an operator expansion of the action.
Moreover, the operator expansion is given by the adjoint action 
so that our ideas concerning invariant integration theory 
apply  \cite{EW}. 
Hilbert space rep\-re\-sen\-ta\-tions of some  cross product algebras 
can be found in \cite{SW} and \cite{EW}.

%\bibliography{matball}
%\bibliographystyle{plain}

\end{document}